\documentclass[final,leqno,onefignum,onetabnum]{siamltexmm}

\usepackage{amssymb}
\usepackage{amsmath}
\DeclareMathAlphabet\gothic{U}{euf}{m}{n}

\usepackage{color}

\def\th{\theta}
\newcommand{\ul}{\textbf}

\newcommand{\R}{\mathbb{R}}

\newcommand{\desda}{\Leftrightarrow}
\newcommand{\SE}{\operatorname{SE(2)}\nolimits}

\newcommand{\eq}[1]{$(\protect\ref{#1})$}

\newtheorem{remark}{Remark}

\title{A PDE Approach to Data-driven Sub-Riemannian Geodesics in SE(2)}

\author{E.J.~Bekkers$^2$, R.~Duits$^{1,2}$, A.~Mashtakov$^2$ and G.R.~Sanguinetti$^1$\thanks{Joint main authors in alphabetic order.
1. Department of Mathematics and Computer Sciences, CASA, TU/e, Eindhoven.
2. Department of Biomedical Engineering,  BMIA, TU/e, Eindhoven}}

\begin{document}
\maketitle
\newcommand{\slugmaster}{%
\slugger{siam}{xxxx}{xx}{x}{x--x}}

\begin{abstract}
We present a new flexible wavefront propagation algorithm for the boundary value problem for sub-Riemannian (SR) geodesics in
the roto-translation group $SE(2) = \mathbb{R}^2 \rtimes S^1$ with a metric tensor depending on a smooth external cost $\mathcal{C}:SE(2) \to [\delta,1]$, $\delta>0$, computed from image data.
The method consists of a first step where a SR-distance map is computed as a viscosity solution of a Hamilton-Jacobi-Bellman (HJB) system derived via Pontryagin's Maximum Principle (PMP).
Subsequent backward integration, again relying on PMP, gives the SR-geodesics. For $\mathcal{C}=1$ we show that our method produces the global minimizers. Comparison with exact solutions shows a remarkable accuracy of the SR-spheres and the SR-geodesics. We present numerical computations of Maxwell points and cusp points, which we again verify for the uniform cost case $\mathcal{C}=1$.
Regarding image analysis applications, tracking of elongated structures in retinal and synthetic images show that our line tracking generically deals with crossings. We show the benefits of including the sub-Riemannian geometry.
\end{abstract}



\pagestyle{myheadings}
\thispagestyle{plain}
\markboth{A PDE approach to Data-Driven Sub-Riemannian geodesics in $SE(2)$}{Bekkers, E.J., Duits, R., Mashtakov, A. and Sanguinetti, G.R.}
\section{Introduction}

In computer vision, it is common to extract salient curves in images via minimal paths or geodesics minimizing a length functional. The minimizing geodesic is defined as the curve that minimizes the length functional, which is typically weighted by a cost function with high values on image locations with high curve saliency. To compute such data-driven geodesics many authors use a two step approach in which firstly a geodesic distance map to a source is computed and then steepest descent on the map gives the geodesics. In a PDE framework, the geodesic map is obtained via wavefront propagation as the viscosity solution of a Hamilton-Jacobi-Bellman equation (the Eikonal equation). For a review of this approach and applications see \cite{sethian,osher,PeyreCohen2010}.

\begin{figure}
\includegraphics[width=\textwidth]{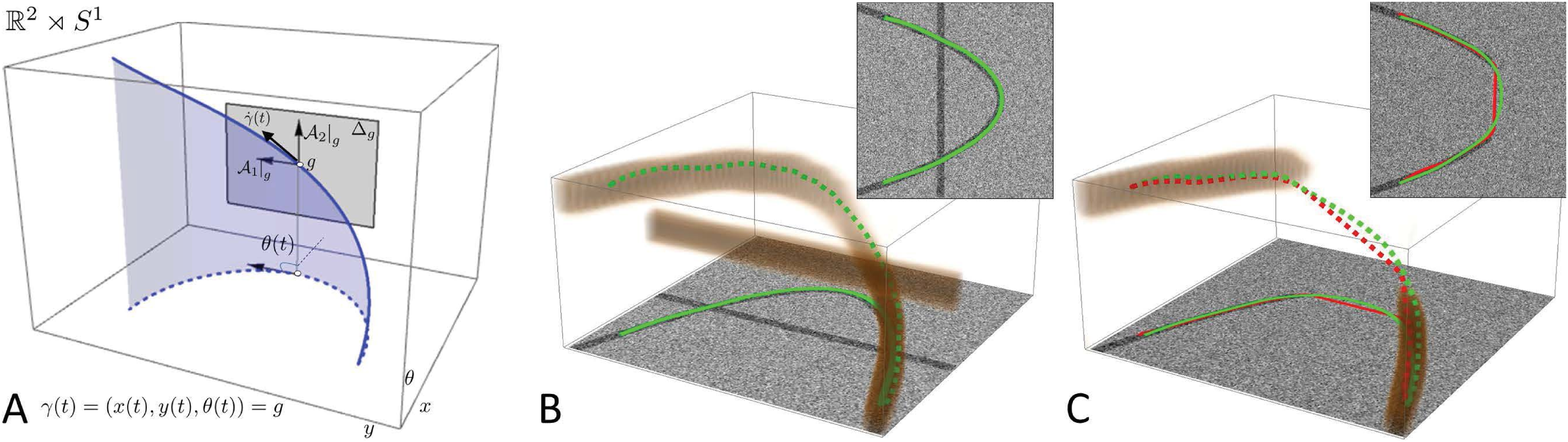}
\centering
\caption{\textbf{A}: Every point in the planar curve $\gamma_{\operatorname{2D}}(t)=(x(t),y(t))$ is lifted to a point $g=\gamma(t)=(x(t),y(t),\theta(t))\in SE(2)$ on an horizontal curve (solid line) by considering the direction of the tangent vector $\dot{\gamma}_{\operatorname{2D}}(t)$ of the planar curve as the third coordinate. Then, tangent vectors $\dot{\gamma}(t) \in \textrm{span}\{\left.\mathcal{A}_{1}\right|_{\gamma(t)},\left.\mathcal{A}_{2}\right|_{\gamma(t)}\} = \left.\Delta\right|_{\gamma(t)}$, Eq.~\!(\ref{leftinvariant}). \textbf{B}: In the lifted domain $SE(2)$ crossing structures are disentangled. \textbf{C}: The SR-geodesic (green) better follows the curvilinear structure along the gap than the Riemannian geodesic (red). }
\label{fig1}
\end{figure}

Another set of geodesic methods, partially inspired by the psychology of vision was developed in \cite{cittisarti,petitot2003}. Here, the roto-translation group $SE(2)=\R^{2} \rtimes S^{1}$ endowed with a sub-Riemannian (SR) metric models the functional architecture of the primary visual cortex and geodesics are used for completion of occluded contours. A stable wavelet-like approach to lift 2D-images to functions on $SE(2)$ was proposed in \cite{DuitsIJCV2007}. Within the $SE(2)$ framework, images and curves are lifted to the 3D space $\mathbb{R}^2 \rtimes S^1$ of coupled positions and orientations in which intersecting curves are disentangled. The SR-structure applies a restriction to so-called horizontal curves which are the curves naturally lifted from the plane (see Fig.~\ref{fig1}A). For explicit formulas of SR-geodesics and optimal synthesis see \cite{yuriSE2FINAL}. SR-geodesics in $SE(2)$ were also studied in \cite{benyosef,boscain,DuitsJMIV,mashtakov,yuriSE2}.
Here, we propose a new wavefront propagation-based method for finding SR-geodesics within $SE(2)$ with a metric tensor depending on a smooth external cost $\mathcal{C}:SE(2) \to [\delta,1]$, $\delta>0$ fixed. Our solution is based on a Hamilton-Jacobi-Bellman (HJB) equation in $SE(2)$ with a SR metric including the cost.
It is of interest to interpret the viscosity solution of the corresponding HJB equation as a sub-Riemannian distance~\cite{trelat}.
Using Pontryagin's Maximum Principle (PMP), we derive the HJB-system
with an Eikonal equation providing the propagation of geodesically equidistant surfaces departing from the origin. We prove this in Thm.~\ref{th:3a}, and we show that SR-geodesics are computed by backtracking via PMP. In Thm.~\ref{th:3}, we consider the uniform cost case (i.e.~$\mathcal{C}=1$) and we show that the surfaces coincide with the SR-spheres, i.e. the surfaces from which every tracked curve is globally optimal. This uniform cost case has been deeply studied in \cite{yuriSE2FINAL} relying on explicit ODE-integration in PMP. In this article, we will rely on a PDE-approach and viscosity solutions of HJB-equations, allowing us to extend to the general case where $\mathcal{C}$ is a smooth cost uniformly bounded from below and above.
We will often use the results in \cite{yuriSE2FINAL} as a golden standard to verify optimality properties of the viscosity solutions and accuracy of the involved numerics of our PDE-approach.
We find a remarkable accuracy compared to exact solutions, 1st Maxwell sets (i.e. the location where for the first time two distinct geodesics with equal length meet), and to the cusp surface \cite{DuitsJMIV,boscain}.

Potential towards applications of the method with non-uniform cost is demonstrated by performing vessel tracking in retinal images. Here the cost function is computed by lifting the images via oriented wavelets, as is explained in Section~\ref{sec:cost}. Similar ideas of computing geodesics via wavefront propagation in the extended image domain of positions and orientations, and/or scales, have been proposed in  \cite{Pechaud,LiYezzi,benmansour}. In addition to these interesting works we propose to rely on a SR geometry.
Let us illustrate some key features of our method. In Fig.~\ref{fig1}B one can see how disentanglement of intersecting structures, due to their difference in orientations, allows to automatically deal with crossings (a similar result can be obtained with the algorithm in \cite{Pechaud}). The additional benefit of using a SR geometry is shown in Fig.~\ref{fig1}C where the SR-geodesic better follows the curvilinear structure along the gap.

\subsection{Structure of the Article}

The article is structured as follows. First, in Section~\ref{sec:Formulation}, we give the mathematical formulation of the curve optimization problem that we aim to solve in this paper. We also provide embedding into previously studied geometrical control problems formulated on $SE(2)$ and/or $\R^{2}$.
In Section~\ref{sec:Solution} we describe our PDE approach that provides the SR geodesic distance map as
the viscosity solution of a boundary value problem (BVP) involving a sub-Riemannian Eikonal equation. Furthermore, in Theorem~\ref{th:3a}, we show that sub-Riemannian geodesics are obtained from this distance map by back-tracing (imposed by the PMP computations presented in Appendix~\ref{app:A}). In Theorem~\ref{th:3} we show that for the uniform cost case (i.e. $\mathcal{C}=1$) such back-tracing will never pass a Maxwell-points or conjugate points, and thereby for $\mathcal{C}=1$ our approach provides only the globally optimal solutions.

In Section~\ref{sec:iterativeproc} we describe an iterative procedure on how to solve the BVP by solving a sequence of initial value problems (IVP) for the corresponding HJB-equation. Before involvement of numerics, we express the exact solutions in concatenated morphological convolutions (erosions) and time-shifts in Appendix~\ref{app:new}. Here we rely on
morphological scale space PDE's \cite{burgeth,akian,DuitsJMIV2013}, and we show that solutions of the iterative procedure converges towards the sub-Riemannian distance map.
Then in Section~\ref{sec:cost} we construct the external cost $\mathcal{C}$, based on a lifting of the original image to an orientation score \cite{DuitsIJCV2007}.

In Section~\ref{ch:implementation}, we describe a numerical  PDE-implementation of our method by using left-invariant finite differences \cite{Franken2009} in combination with an upwind-scheme \cite{Rouy}. Finally, in Section~\ref{sec:Experiments}, we present a verification of the proposed method by comparison with the exact solutions. We also provide simple numerical approaches (extendable to the non-uniform cost case) to compute 1st Maxwell points, conjugate points and cusp points \cite{DuitsJMIV}, which we verify for the uniform cost case with results in \cite{yuriSE2FINAL}. Finally, in Subsection~\ref{ch:appl}, we show application of the method to vessel tracking in optical images of the retina.
We discuss the two main parameters that are involved: the balance between external and internal cost, and the balance between spatial and angular motion. First feasibility studies are presented on patches, and we discuss on how to proceed towards automated retinal vessel tree segmentation.

\section{Problem Formulation}\label{sec:Formulation}

The roto-translation group $SE(2)$ carries group product:
\[
g g'=(\ul{x},R_{\theta}) (\ul{x}',R_{\theta'})=
(R_{\theta}\ul{x}' +\ul{x},R_{\theta+\theta'}).
\]
where $R_{\theta}$ is a counter-clockwise planar rotation over angle $\theta$.
This group can be naturally identified with the
coupled space of positions and orientations $\R^{2} \rtimes S^{1}$, by identifying
$R_{\theta} \leftrightarrow \theta$ while imposing  $2\pi$-periodicity on $\theta$.
Then for each $g \in SE(2)$ we have the left multiplication $L_{g}h=gh$. Via the push-forward $(L_{g})_*$
of the left-multiplication we get the left-invariant vector fields $\{\mathcal{A}_{1},\mathcal{A}_{2},\mathcal{A}_{3}\}$ from the Lie-algebra basis $\{A_{1},A_{2},A_{3}\}=\{\left.\partial_{x}\right|_{e}, \left.\partial_{\theta}\right|_{e},
\left. \partial_{y} \right|_{e}\}$
at the unity $e=(0,0,0)$:
\begin{equation} \label{leftinvariant}
\begin{array}{l}
 \left.\mathcal{A}_{1}\right|_{g}= \cos \theta \left.\partial_{x}\right|_{g} +\sin \theta
 \left.\partial_{y}\right|_{g}= (L_{g})_* \left.\partial_{x}\right|_{e}, \\
 \left.\mathcal{A}_{2}\right|_{g}= \left.\partial_{\theta}\right|_{g}=(L_{g})_* \left.\partial_{\theta}\right|_{e}, \\
 \left. \mathcal{A}_{3}\right|_{g}= -\sin \theta \left.\partial_{x}\right|_{g} +\cos \theta
 \left.\partial_{y}\right|_{g}= (L_{g})_* \left.\partial_{y}\right|_{e}.
\end{array}
\end{equation}
Then all tangents $\dot{\gamma}(t) \in T_{\gamma(t)}(SE(2))$ along smooth curves
$t\mapsto \gamma(t)=(x(t),y(t),\theta(t)) \in SE(2)$ can be expanded as
$\dot{\gamma}(t)=\sum_{k=1}^{3} u^{k}(t) \left.\mathcal{A}_{k}\right|_{\gamma(t)}$,
where the contravariant components $u^k(t)$ of the tangents (velocities) can be considered as the control variables.
Not all curves $t \mapsto \gamma(t)$ in $SE(2)$ are naturally lifted from the plane in the sense that $\theta(t)=\arg(\dot{x}(t)+i\, \dot{y}(t))$.
This holds for so-called \emph{horizontal curves} which have $u^3=0$ and thus
$\dot{\gamma}(t)=\sum_{k=1}^{2} u^{k}(t) \left.\mathcal{A}_{k}\right|_{\gamma(t)}$.
The allowed (horizontal) directions in tangent bundle $T(SE(2))$ (see Fig.~\ref{fig1}A) form a so-called distribution
\[
\Delta:=\textrm{span}\{\mathcal{A}_{1},\mathcal{A}_{2}\}.
\]
Therefore we consider SR-manifold~\cite{mont} $(SE(2),\Delta,G^{\mathcal{C}})$, with \\
$G^{\mathcal{C}}: SE(2)\times \Delta \times \Delta \to \R$ denoting the inner product given by
\begin{equation} \label{metric}
G^{\mathcal{C}}|_{\gamma(t)}(\dot{\gamma}(t),\dot{\gamma}(t)) = \mathcal{C}^2\left(\gamma(t)\right)\left(\beta^2 |\dot{x}(t)\cos\th(t)  \!+\!\dot{y}(t)\sin\th(t) |^2 + |\dot{\th}(t)|^2\right),
\end{equation}
with $\gamma:\R\! \to\! SE(2)$ a smooth curve on $\R^{2}\rtimes S^{1}$, $\beta>0$ constant, $\mathcal{C}:SE(2)\to [\delta,1]$ the \emph{given external smooth cost which is bounded from below by $\delta>0$}. 
\begin{remark}
Define $\mathcal{L}_{g}\phi(h)=\phi(g^{-1}h)$ then we have:
$$
G^{\mathcal{C}}|_{\gamma}(\dot{\gamma},\dot{\gamma})= \left.G^{\mathcal{L}_{g}\mathcal{C}}\right|_{g \gamma}\left( \;(L_{g})_{*}\dot{\gamma}, (L_{g})_{*}\dot{\gamma}\;\right).
$$
Thus, $G^{\mathcal{C}}$ is not left-invariant, but if shifting the cost as well, we can, for the computation of SR-geodesics, restrict ourselves to $\gamma(0)=e$.
\end{remark}

We study the problem of finding SR minimizers, i.e. for given boundary conditions
$\gamma(0) = e, \gamma(T) = g$,
we aim to find the horizontal curve $\gamma(t)$ (having $\dot{\gamma} \in \Delta$) that minimizes the total SR length
\begin{equation}
l = \int_0^{T} \sqrt{\left.G^{\mathcal{C}}\right|_{\gamma(t)}(\dot{\gamma}(t),\dot{\gamma}(t))}\, {\rm d} t.
\end{equation}
If $t$ is the SR arclength parameter, our default parameter,
then $\sqrt{\left. G^{\mathcal{C}}\right|_{\gamma(t)}(\dot{\gamma}(t),\dot{\gamma}(t))}=1$ and $l=T$.
Then, SR minimizers $\gamma$ are solutions to the optimal control problem (with free $T>0$):
\begin{equation} \label{problemP}
\textbf{P}_{\textbf{mec}}^{\mathcal{C}}(SE(2))\;:\; \left\{
\begin{array}{c}
\dot{\gamma} = u^1 \left.\mathcal{A}_1\right|_{\gamma} + u^2 \left. \mathcal{A}_{2}\right|_{\gamma},  \\
\gamma(0) = e, \qquad \gamma(T) = g, \\
l(\gamma(\cdot)) = \int_{0}^{T} \mathcal{C}(\gamma(t))\sqrt{\beta^2 |u^1(t)|^2 +  |u^2(t)|^2}\, {\rm d}t \to \min, \\
\qquad \gamma(t) \in \SE, \quad (u^1(t), u^2(t)) \in \R^2, \quad \beta>0.
\end{array}
\right.
\end{equation}
In the naming of this geometric control problem we adhere to terminology in previous work \cite{boscain,DuitsJMIV}.
Stationary curves of the problem (\ref{problemP}) are found via PMP~\cite{notes}. Existence of minimizers follows from Chow-Rashevsky and Filippov's theorem~\cite{notes}, and because of absence of abnormal trajectories
(due to the 2-bracket generating distribution $\Delta$) they are smooth.
\begin{remark} \label{rem:2}
The Cauchy-Schwarz inequality implies that the minimization problem for the SR length functional $l$ is equivalent (see e.g.~\cite{mont}) to the minimization problem for the action functional with fixed $T$:
\begin{equation}\label{action}
J(\gamma) = \frac12 \int_0^T \mathcal{C}^2(\gamma(t))(\beta^2 |u^1(t)|^2 +|u^2(t)|^2)\ {\rm d}t.
\end{equation}
\end{remark}

\subsection{Embedding into Geometric Control Theory}

The problem (\ref{problemP}) actually comes from a mechanical problem in geometric control, where a so-called Reeds-Shepp car \cite{reeds} can proceed both forward and backward in the path-optimization. As pointed out in \cite{charlot} such a problem, for certain end conditions, cannot be considered as a curve optimization problem on the plane. In fact in \cite{DuitsJMIV2013} the set $\gothic{R} \subset SE(2)$ of allowable end-conditions $g$ is explicitly determined for $\mathcal{C}=1$. The boundary of the set $\gothic{R}$ is formed by end-points of geodesics whose spatial projections either depart or end in a cusp
\cite[Fig.12,Thm.9]{DuitsJMIV2013}. For $g=(\ul{x},\theta) \in \gothic{R} \subset SE(2)$ and $\mathcal{C}=1$ the following problem is well-posed
\begin{equation} \label{problemPcurve}
\textbf{P}_{\textbf{curve}}^{\mathcal{C}}(\R^{2})\;:\;
\left\{
\begin{array}{l}
\gamma(0) = \ul{0}, \  \gamma(L) = \ul{x}, \\
\dot\gamma(0) = (1,0)^T, \ \dot{\gamma}(L)=(\cos \theta,\sin \theta)^T, \\
l(\gamma(\cdot)) = \int \limits_{0}^{L} \mathcal{C}(\gamma(s))\, \sqrt{\beta^2 + \kappa^{2}(s)}\, {\rm d}s \to \min, \\
\gamma:[0,L] \to \R^{2},  \beta>0,
\end{array}
\right.
\end{equation}
where $L$ denotes spatial length and $\kappa$ curvature of smooth curve $\gamma \in C^{\infty}([0,L],\R^{2})$.
For $\mathcal{C} \neq 1$ the set of end-conditions $(\ul{x},\theta)$ for which the above problem is well-posed
varies with $\mathcal{C}$.

\begin{figure}
\includegraphics[width=\hsize]{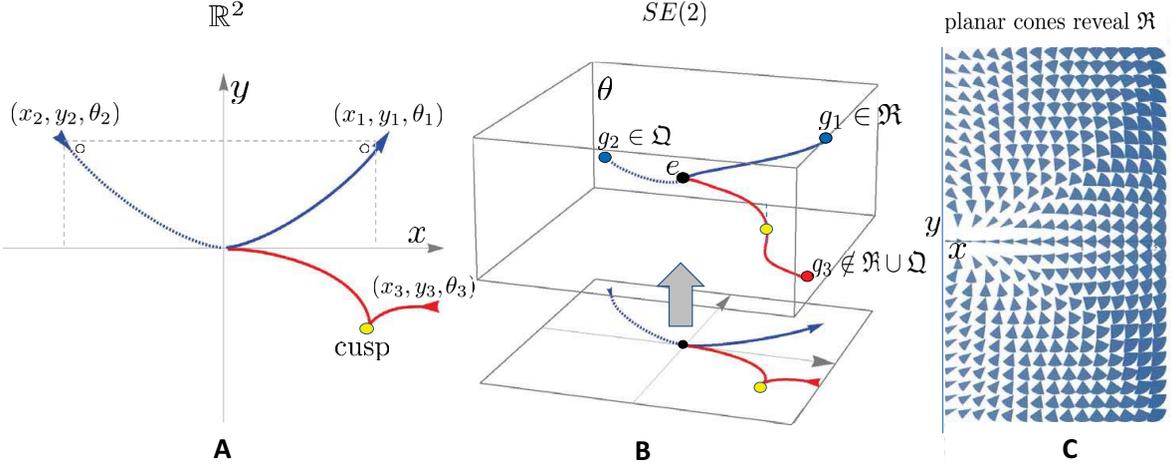}
\caption{3 Sub-Riemannian geodesics for uniform cost $\mathcal{C}=1$ to reveal the differences of the 3 geometric control problems (\ref{problemP}), (\ref{problemPcurve}) and (\ref{problemPcontour}), for $\beta=1$.
\textbf{A}: plots of spatial projections of sub-Riemannian geodesics in $\R^2$,
\textbf{B}: SR-geodesics in $SE(2)$, \textbf{C}: A part of $\gothic{R}$ (the set of end-conditions for which $\textbf{P}_{\textbf{curve}}(\R^2)$ is well-defined \cite{DuitsJMIV}) depicted as reachable cones around the origin.
End condition $g_{1}=(x_1,y_1,\theta_1) \in \gothic{R}$ yields the minimizing curve of $\textbf{P}_{\textbf{contour}}(SE(2))$, $\textbf{P}_{\textbf{mec}}(SE(2))$
and $\textbf{P}_{\textbf{curve}}(\R^2)$. End condition $g_{2}=(-x_1,y_1,-\theta_{1})$ yields the minimizing curve in $\textbf{P}_{\textbf{contour}}(SE(2))$, $\textbf{P}_{\textbf{mec}}(SE(2))$, and it is invalid for $\textbf{P}_{\textbf{curve}}(\R^2)$. End condition $g_3$ is invalid for both $\textbf{P}_{\textbf{curve}}(\R^2)$ and
$\textbf{P}_{\textbf{contour}}(SE(2))$, as it induces an internal cusp.
For $\mathcal{C}=1$ the set of allowable end conditions for $\textbf{P}_{\textbf{contour}}(SE(2))$ equals $\gothic{R}^{\mathcal{C}=1}=
\gothic{R} \cup \gothic{Q}$. For $\mathcal{C}\neq 1$ this set $\gothic{R}^{\mathcal{C}}$ differs, and can be computed.
\label{fig:Ps}}
\end{figure}

In this article we take $\textbf{P}_{\textbf{mec}}^{\mathcal{C}}(SE(2))$ (and not (\ref{problemPcurve})) as a venture point, as for the detection of elongated structures we must look both forward and backward.
From the image analysis point of view, our primary interest is in the following refined problem of $\textbf{P}_{\textbf{mec}}^{\mathcal{C}}(SE(2))$:
\begin{equation} \label{problemPcontour}
\textbf{P}_{\textbf{contour}}^{\mathcal{C}}(SE(2))\;:\; \left\{
\begin{array}{l}
\dot{\gamma}(t) = u^1(t) \left.\mathcal{A}_1\right|_{\gamma(t)} + u^2(t) \left. \mathcal{A}_{2}\right|_{\gamma(t)}, \textrm{ for }t\in [0,T]  \\[6pt]
\gamma(0) = e, \qquad \gamma(T) = g =(\ul{x},\theta) \in \gothic{R}^{\mathcal{C}}, \\
l(\gamma(\cdot)) = \int \limits_{0}^{T} \mathcal{C}(\gamma(t))\sqrt{\beta^2 |u^1(t)|^2 +  |u^2(t)|^2}\, {\rm d}t \to \min, \\
\textrm{with curve }\gamma: [0,T]\to \SE, \textrm{ with controls: } \\ (u^1(t), u^2(t)) \in \R^2, \textrm{ and }u^{1}(t)\textrm{ does not change sign},
\end{array}
\right.
\end{equation}
where $\gothic{R}^{\mathcal{C}}$ is the set of all $g \in SE(2)$ such that the minimizing SR-geodesic(s) $\gamma(\cdot)=(\ul{x}(\cdot),\theta(\cdot))$ do not exhibit a cusp in their spatial projections $\ul{x}(\cdot)$.
\begin{remark}
For existence results of the minimizers, a priori $\mathbb{L}_{1}$-conditions have to be imposed on the controls $u_{1}$, $u_{2}$. However, in retrospect, after application of a generalized version of PMP \cite{Vinter}, \cite[Thm 5.2 and App.B]{boscain}, the controls of the minimizers turn out to be smooth and bounded. Regarding extension to $\mathcal{C} \neq 1$ (recall (\ref{metric})) we note that $\delta \leq \mathcal{C} \leq 1$.
\end{remark}
\begin{remark}
For $\mathcal{C}=1$ one has
\[
\gothic{R}^{\mathcal{C}=1}= \gothic{R} \cup \gothic{Q} \textrm{ with }
\gothic{Q}=\{(x,y,\theta) \in SE(2)\;|\; (-x,y,-\theta) \in \gothic{R}\},
\]
and with $\gothic{R} \subset \{(x,y,\theta) \in SE(2)\;|\; x \geq 0\}$ described explicitly in \cite{DuitsJMIV}, and partly depicted in item C of Figure~\ref{fig:Ps}.
In Section~\ref{ch:expcusp} 
we shall provide a very simple numerical tool to compute the surface in $SE(2)$ where cusps appear also for $\mathcal{C}\neq 1$.
This surface is a boundary of a volume in $SE(2)$ that contains the set
$\gothic{R}^{\mathcal{C}}$.
\end{remark}
\ \\
Summarizing, in $\textbf{P}_{\textbf{mec}}^{\mathcal{C}}(SE(2))$ the optimal control $t \mapsto u^{1}(t)$ can switch sign.
\begin{itemize}
\item If $g$ is chosen such that the optimal control $u^{1} \geq 0$ then the lift of problem $\textbf{P}_{\textbf{curve}}^{\mathcal{C}}(\R^2)$ coincides with
$\textbf{P}_{\textbf{mec}}^{\mathcal{C}}(SE(2))$ and also with $\textbf{P}_{\textbf{contour}}^{\mathcal{C}}(SE(2))$.
\item If $g$ is chosen such that the optimal control $u^{1} \leq 0$ then problem $\textbf{P}_{\textbf{mec}}^{\mathcal{C}}(SE(2))$ and problem $\textbf{P}_{\textbf{contour}}^{\mathcal{C}}(SE(2))$ coincide.
\item If $g$ is chosen such that the optimal control $u^{1}(t)$ switches sign at some internal time $t \in (0,T)$, then $g \in SE(2)\setminus \gothic{R}^{\mathcal{C}}$ and the spatial projection of the corresponding minimizing SR-geodesic(s) has an internal cusp, which we consider not desirable in our applications of interest.
\end{itemize}
See Figure~\ref{fig:Ps}, where each of the 3 above cases is illustrated.

Now let us return to our main objective, that is to derive solutions of (\ref{problemP}) via a PDE-wavefront propagation method for data-driven SR-geodesics.
The word `data-driven' refers to the fact that our PDE-approach is also applicable to $\mathcal{C} \neq 1$, a property that does not apply to previous methodology following ODE-approaches \cite{yuriSE2FINAL,boscain,mashtakov,DuitsJMIV,petitot2003} on SR-geodesics in vision/image analysis.

\section{Solutions via Data-driven Wavefront Propagation}\label{sec:Solution}

\begin{figure}
\centerline{
\includegraphics[width=0.8\textwidth]{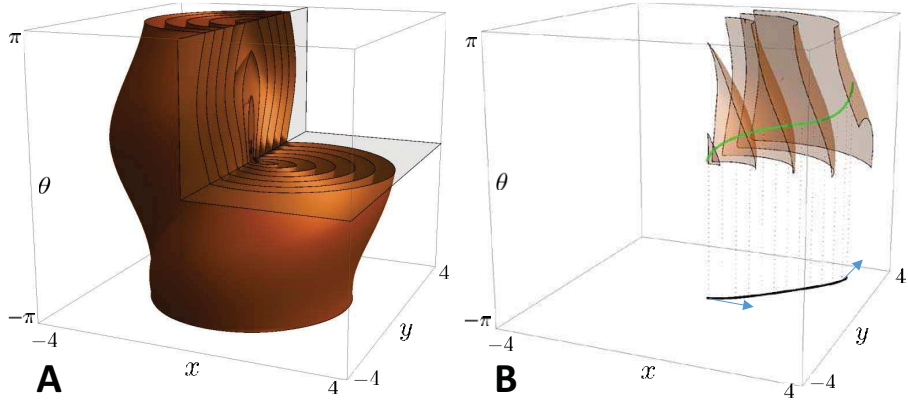}
}
\caption{
\textbf{A-B}: Our method provides both geodesically equidistant surfaces $\mathcal{S}_t$ (\ref{St}) (depicted in A) and SR-geodesics. As depicted in B: geodesic equidistance holds with unit speed for all SR-geodesics passing through the surface, see Thm~\ref{th:3a}.
Via Thm.~\ref{th:3} we have that $W(g)=d(g,e)$ and $\{\mathcal{S}_{t}\}_{t \geq 0}$ is the family of SR-spheres
with radius $t$ depicted in this figure for the uniform cost case. 
}
\label{figWaveFronts}
\end{figure}

The following theorem summarizes our method for the computation of data-driven sub-Riemannian geodesics in $SE(2)$. The idea is illustrated in Fig.~\ref{figWaveFronts}.
\begin{theorem} \label{th:3a}
Let $W(g)$ be a solution of the following boundary value problem (BVP) with Eikonal-equation
\begin{equation} \label{Eik}
\left\{
\begin{array}{l}
\sqrt{(\mathcal{C}(g))^{-2}
\left(\beta^{-2}|\mathcal{A}_{1}W(g)|^2 + |\mathcal{A}_{2}W(g)|^2\right)}-1=0, \textrm{ for }g\neq e, \\
W(e)=0.
\end{array}
\right.
\end{equation}
Then the iso-contours
\begin{equation} \label{St}
\mathcal{S}_{t}=\{g \in SE(2)\;|\; W(g)=t\}
\end{equation}
are geodesically equidistant with speed
$\frac{dt}{dt}=\mathcal{C}(\gamma(t)) \sqrt{\beta^2 |u^{1}(t)|^2 + |u^{2}(t)|^2}=1$ and they
provide a specific part of the SR-wavefronts departing from
$e=(0,0,0)$. A
SR-geodesic departing from $g\in SE(2)$ is found by backward integration
\begin{equation}\label{steepest}
\dot{\gamma}_b(t)= - \frac{\left.\mathcal{A}_{1}W\right|_{\gamma_b(t)}}{(\beta \, \mathcal{C}(\gamma_b(t)))^2} \left.\mathcal{A}_{1}\right|_{\gamma_b(t)}-
\frac{\left.\mathcal{A}_{2}W\right|_{\gamma_b(t)}}{
(\mathcal{C}(\gamma_b(t)))^2} \left.\mathcal{A}_{2}\right|_{\gamma_b(t)}, \;\;\;\;\;\;\gamma_{b}(0)=g.
\end{equation}
\end{theorem}
\textbf{Proof }
The definition of geodesically equidistant surfaces is given in Definition~\ref{def:geodS} in Appendix~\ref{app:C}.
Furthermore, in Appendix~\ref{app:C} we provide two lemmas needed for the proof.
In Lemma~\ref{th:1}, we
connect the Fenchel transform on $\Delta$, to the Fenchel transform on $\R^{2}$
to obtain the result on geodesically equidistant surfaces in $(SE(2),\Delta,G^{\mathcal{C}})$.
Then, in Lemma~\ref{th:2} in Appendix~\ref{app:C}, we derive the HJB-equation for the homogeneous Lagrangian as a limit from the HJB-equation for the squared Lagrangian. The back-tracking result follows
from application of PMP to the equivalent action functional formulation (\ref{action}). Akin to the $\R^{d}$-case~\cite{bressan}, characteristics are found by integrating the ODE's of the PMP where according to the proof of Lemma~\ref{th:1} we must set $p={\rm d}^{SR}W$, see Remark~\ref{rem:3} below.
$\hfill \Box$

The next theorem provides our main theoretical result. Recall that Maxwell points are $SE(2)$ points where two distinct geodesics with the same length meet.
The 1st Maxwell set corresponds to the set of Maxwell-points where the distinct geodesics meet for the first time.
In the subsequent theorem we will consider a specific solution to (\ref{Eik}), namely the viscosity solution
as defined in Definition~\ref{def:VSEik} in Appendix~\ref{app:B}.
\begin{theorem} \label{th:3}
Let $\mathcal{C}=1$.
Let $W(g)$ be the viscosity solution of the BVP (\ref{Eik}). Then $\mathcal{S}_{t}$, Eq.~\!(\ref{St}), equals the
SR-sphere of radius $t$.
Backward integration via (\ref{steepest})
provides globally optimal geodesics reaching $e$
at $t=W(g)=d(g,e):=$
\begin{equation} \label{SRdist}
\min_{{\small \begin{array}{l}
\gamma \in C^{\infty}(\mathbb{R}^+,SE(2)), T\geq 0, \\
\dot{\gamma} \in \Delta, \gamma(0)=e,\gamma(T)=g
\end{array}}}\!
\int_{0}^{T}\sqrt{|\dot{\theta}(\tau)|^2 +\beta^{2}
|\dot{x}(\tau)\cos \theta(\tau)\!+\!\dot{y}(\tau)\sin \theta(\tau)|^2}\; {\rm d}\tau,
\end{equation}
and $\gamma_b(t)=\gamma^{min}(d(g,e)-t)$.
The SR-spheres $\mathcal{S}_{t}=\{g \in SE(2) \;|\; d(g,e)=t\}$ are non-smooth 
at the 1st Maxwell set $\mathcal{M}$, cf.\!~\cite{yuriSE2FINAL}, 
contained in
\begin{equation} \label{Maxwell}
\begin{array}{l}
\mathcal{M}\subset\left\{(x,y,\theta) \in SE(2)\;|\; x\cos \frac{\theta}{2}+
y \sin \frac{\theta}{2}=0 \vee  \; \theta=\pi   \right\},
\end{array}
\end{equation}
and the back-tracking (\ref{steepest}) does not pass the 1st Maxwell set.
\end{theorem}

\begin{figure}
\centerline{
\includegraphics[width=\hsize]{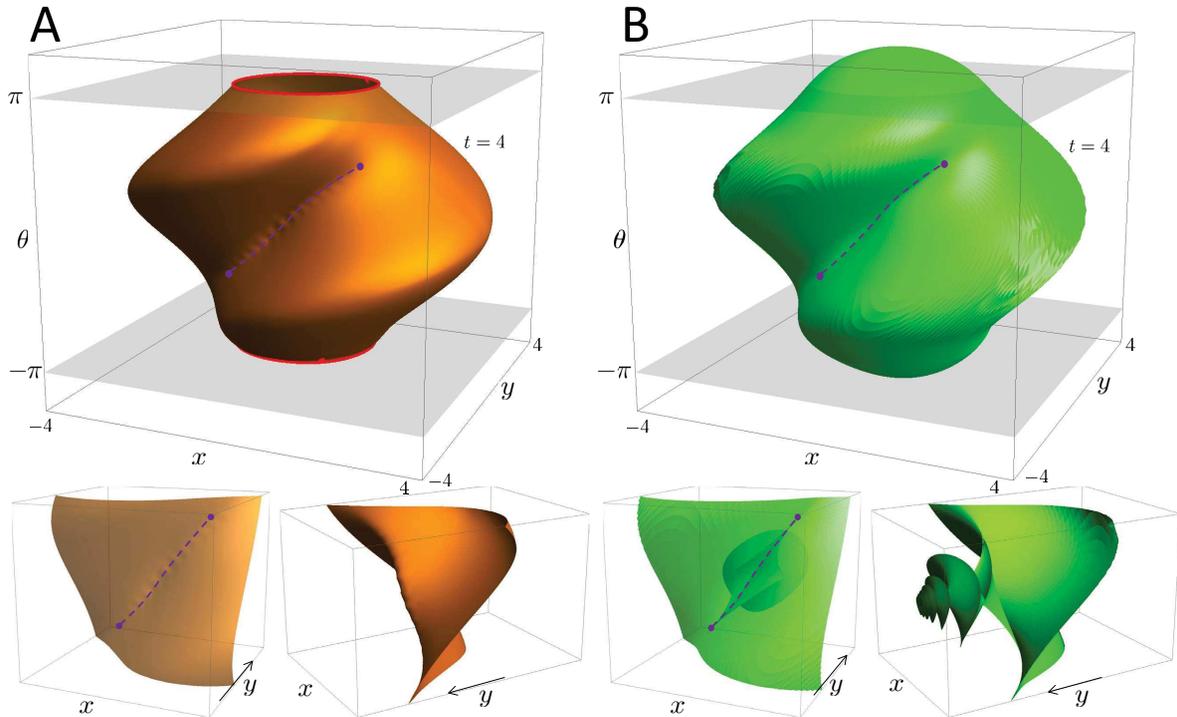}
}
\caption{\textbf{A: }SR-sphere $\mathcal{S}_{t}$ for $t=4$ obtained by the method
in Thm.~\ref{th:3a} using $\mathcal{C}=1$ and $\delta^{M}_{e}$ as initial condition via viscosity solutions
of the HJB-equation (\ref{thePDE}) implemented according to Section~\ref{ch:implementation}. \textbf{B: }The full SR-wavefront departing from $e$ via the method of
characteristics and formulae in \cite{yuriSE2} giving rise to interior folds
(corresponding to multiple valued non-viscosity solutions of the HJB-equation). The Maxwell set $\mathcal{M}$ consists precisely
of the dashed line on $x\cos \frac{\theta}{2} + y \sin \frac{\theta}{2}=0$ and the red circles at $|\theta|=\pi$.
The dots are 2 (of the 4) conjugate points on $\mathcal{S}_{t}$
which are limits of 1st Maxwell points (but not Maxwell points themselves). In \textbf{B} we see the astroidal
structure of the conjugate locus \cite{yuriSE2CONJ,chakir}. In \textbf{A} we see that the
\emph{unique} viscosity solutions stop
at the 1st Maxwell set. Comparison of \textbf{A} and \textbf{B} shows the global optimality and accuracy of our method at
\textbf{A}.
\label{fig:gonz}}
\end{figure}

\begin{figure}
\centerline{
\includegraphics[width=\textwidth]{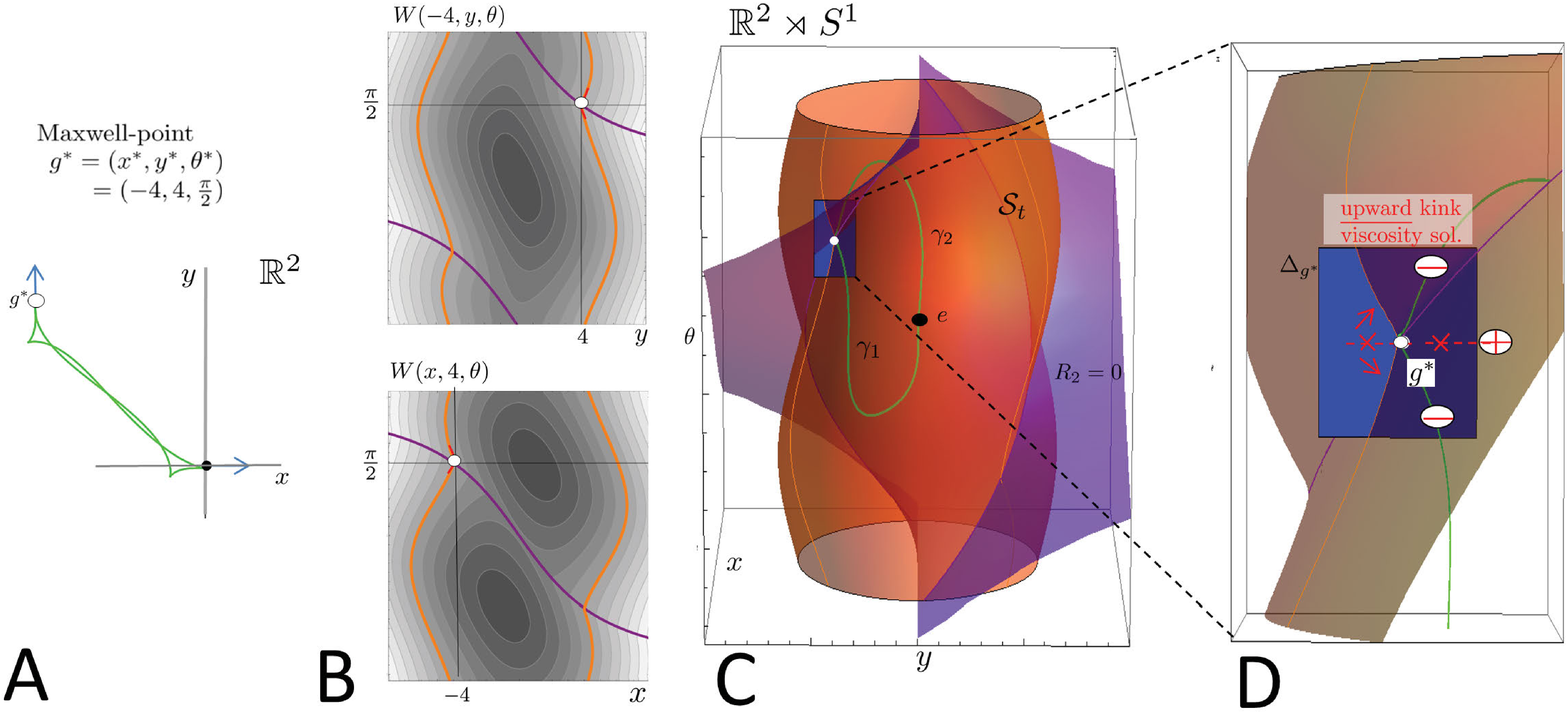}
}
\caption{Maxwell point $g^{*}=(-4,4,\pi/2)$ (in white) on SR-sphere
$\mathcal{S}_{t}$ (in orange) for $\mathcal{C}=1$.
At $g^*$ two SR-geodesics $\gamma_{1} \not\equiv \gamma_{2}$ with equal SR-length $t$ meet ($\gamma_1(t)=\gamma_2(t)$).
From left to right: \textbf{A}: projection of $\gamma_1$ and $\gamma_2$ on the plane $(x,y)$, \textbf{B}:
2D-slices ($x=x^*$, $y=y^*$) of level sets of $W(g)$ with distinguished value $W(g)=t$ (again in orange).
On top we plotted, the Maxwell point, the intersection of surface $x \cos \frac{\theta}{2} + y \sin \frac{\theta}{2}=0$ (in purple,
this set contains a part of the 1st Maxwell set)  with the 2D-slices.
\textbf{C}: The SR-sphere $\mathcal{S}_t$ in $SE(2)$, \textbf{D}: section around
$g^{*}$ revealing the upward kink due to the viscosity solution.
From this kink we see that the tracking (\ref{steepest}) does not cross a 1st Maxwell point as indicated in
red, yielding global optimality in Thm.~\ref{th:3}.}
\label{globaloptim}
\end{figure}

Proof of Thm.~\ref{th:3} can be found in Appendix~\ref{app:D}. The global optimality and non-passing of the 1st Maxwell set can be observed in Fig.~\ref{fig:gonz}. For the geometrical idea of the proof see Fig.~\ref{globaloptim}.
\begin{remark} \label{rem:3}
The stationary solutions of (\ref{thePDE}) satisfy the SR-Eikonal equation (\ref{Eik}).
The Hamiltonian $H^{fixed}$ for the equivalent fixed time problem (\ref{action})
equals
\begin{equation}\label{Hfixed}
H^{fixed}(g,p)= (\mathcal{C}(g))^{-2}\left(\beta^{-2}h_{1}^2 + h_2^2\right)=1/2,
\end{equation}
with momentum covector $p=h_1\omega^{1} +h_{2}\omega^2 +h_{3} \omega^3$ expressed in dual basis $\{\omega^{i}\}_{i=1}^{3}$
given by
\begin{equation} \label{dualbasis}
\langle \omega^{i},\mathcal{A}_{j}\rangle=\delta^{i}_j \desda \omega^{1}=\cos \theta {\rm d}x +\sin \theta {\rm d}y,\  \omega^{2}= {\rm d}\theta,\ \omega^{3}=-\sin \theta {\rm d}x+\cos \theta {\rm d}y.
\end{equation}
The Hamiltonian $H^{free}$ for the free time problem (\ref{problemP}) minimizing $l$ equals
\begin{equation} \label{Hfree}
H^{free}(g,p)=\sqrt{2H^{fixed}(g,p)}-1=0.
\end{equation}
For details see Appendix \ref{app:A} and \ref{app:C}.
Eq.~\!(\ref{Eik}) can be written as $H^{free}(g,p)=0$ with momentum\footnote{Note that the sub-Riemannian gradient $\nabla^{SR}W= G^{-1} {\rm d}W=\sum_{i=1}^{2}\beta_{i}^{-2}\mathcal{A}_{i}W \mathcal{A}_{i}$, with $\beta_{1}=\beta$, $\beta_{2}=1$, by definition is the Riesz-representative (being a vector) of this SR-derivative (being a covector).} 
\[
p = {\rm d}^{SR}W:=\sum_{i=1}^{2}(\mathcal{A}_i W)\;\omega^i.
\]
\end{remark}
\begin{remark} \label{rem:MW}
SR geodesics loose their optimality either at a Maxwell point or at a conjugate point
(where the integrator of the canonical ODE's, mapping initial momentum $p_0$ and time $t>0$ to end-point $\gamma(t)$, is degenerate \cite{notes}). Some conjugate points are limits of Maxwell points, see Fig.~\ref{fig:gonz}, where the 1st astroidal shaped conjugate locus coincides with the void regions (cf.~\!\cite[fig.1]{void}) after 1st 
Maxwell set $\mathcal{M}$. When setting a Maxwell point as initial condition, the initial derivative $\left. {\rm d}^{SR}W \right|_{\gamma_{b}(0)}$ is not defined. Here there are 2 horizontal directions with minimal slope, taking these directions our algorithm produces the results in Fig.~\ref{globaloptim}A and Fig.~\ref{fig:MaxwellDemo}.
\end{remark}

\section{An Iterative IVP-procedure to Solve the SR-Eikonal BVP}\label{sec:iterativeproc}

To obtain an iterative implementation to obtain the viscosity solution of the SR-Eikonal BVP given by\!~(\ref{Eik}), we rely on viscosity solutions of initial value problems (IVP).
In this approach we put a connection between morphological scale spaces \cite{burgeth,akian}, and morphological convolutions with morphological kernels, on the SR manifold $(SE(2),\Delta,G^{\mathcal{C}})$ and the SR Eikonal BVP.

In order to solve the sub-Riemannian Eikonal BVP~(\ref{Eik}) we resort to subsequent auxiliary IVP's on $SE(2)$ for each $r \in [r_n,r_{n+1}]$, with $r_n= n \epsilon$ at step $n \in \mathbb{N} \cup\{0\}$, $\epsilon>0$ fixed:
\begin{equation} \label{thePDE}
\left\{
\begin{array}{ll}
\frac{\partial W^{\epsilon}_{n+1}}{\partial r}(g,r) &= 1-\sqrt{(\mathcal{C}(g))^{-2}
\left(\beta^{-2}|\mathcal{A}_{1}W^{\epsilon}_{n+1}(g,r)|^2 + |\mathcal{A}_{2}W^{\epsilon}_{n+1}(g,r)|^2\right)},
 \\
W^{\epsilon}_{n+1}(g,r_n)&= W^{\epsilon}_{n}(g,r_n) \textrm{ for }g \neq e, \\
W^{\epsilon}_{n+1}(e,r_n)&=0
\end{array}
\right.
\end{equation}
for $n=1,2,\ldots$, and
\begin{equation} \label{thePDE1}
\left\{
\begin{array}{ll}
\frac{\partial W^{\epsilon}_{1}}{\partial r}(g,r) &= 1-\sqrt{(\mathcal{C}(g))^{-2}
\left(\beta^{-2}|\mathcal{A}_{1}W^{\epsilon}_{1}(g,r)|^2 + |\mathcal{A}_{2}W^{\epsilon}_{1}(g,r)|^2\right)},
 \\
W^{\epsilon}_{1}(g,0)&= \delta_{e}^{M}(g),
\end{array}
\right.
\end{equation}
for $n=0$, where $\delta_{e}^{M}$ is the morphological delta given by
\[
\delta_{e}^{M}(g) =
\left\{
\begin{array}{ll}
0 & \textrm{ if }g=e \\
\infty &\textrm{else. }
\end{array}
\right.
\]
Let $W_{n+1}^{\epsilon}$ denote the viscosity solution of (\ref{thePDE}) carrying the following support
\[
\textrm{supp}(W_{n+1}^{\epsilon})= SE(2) \times [r_n,r_{n+1}], \textrm{ with }r_n=n\epsilon,
\]
so in (\ref{thePDE}) at the n-th iteration ($n\geq 1$) we use, for $g\neq e$, the end condition $W_{n}^{\epsilon}(g,r_n)$ of the $n$-th evolution as an initial condition $W_{n+1}^{\epsilon}(g,r_n)$ of the $(n+1)$-th evolution. Only for  $g=e$ we set initial condition $W_{n+1}^{\epsilon}(e,r_n)=0$. Then
we define the pointwise limit
\begin{equation} \label{stationary}
W^{\infty}(g):=\lim \limits_{\epsilon \to 0}\left( \lim \limits_{n \to \infty} W^{\epsilon}_{n+1}(g,n\epsilon)\right).
\end{equation}
Finally, regarding the application of our optimality results, it is important that each IVP-solution $W^{\epsilon}_{n+1}(g,r)$ is the unique \emph{viscosity} solution of (\ref{thePDE}), as then via (\ref{stationary}) the viscosity property for the viscosity solutions of the HJB-IVP problem naturally carries over to the viscosity property of the viscosity solutions of system (\ref{Eik}). Thus we obtain $W=W^{\infty}$ as the unique viscosity solution of the
SR-Eikonal BVP.

\begin{figure}
\includegraphics[width=\hsize]{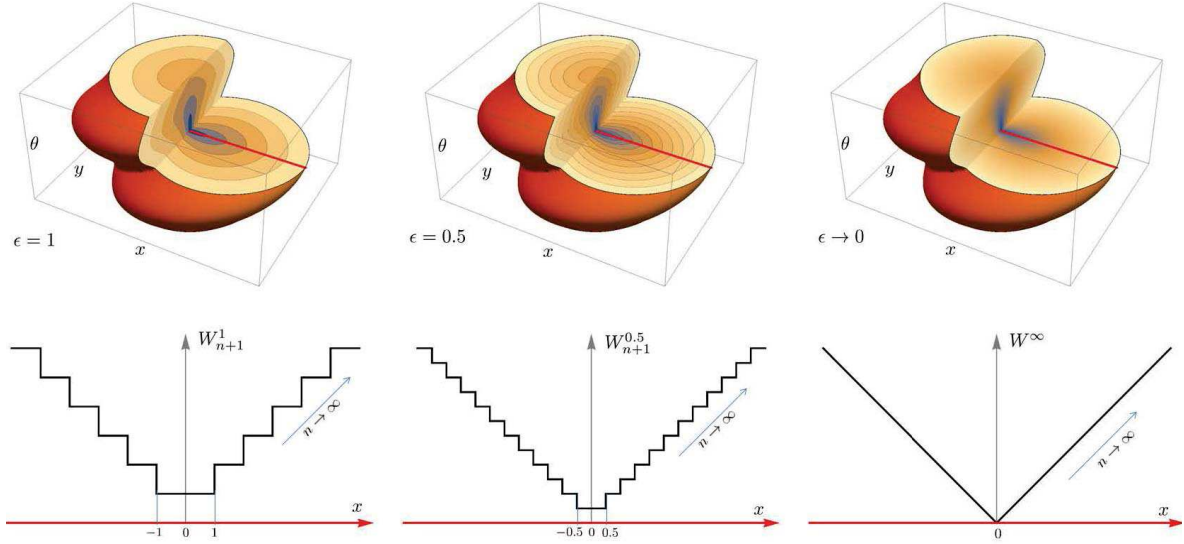}
\caption{Illustration of the pointwise limits in Eq.~\!(\ref{stationary}).
Top:
plot of $g \mapsto \lim \limits_{n \to \infty} W_{n+1}^{\epsilon}(g,r_{n+1})$ (from left to right, resp. for $\epsilon=1$, $\epsilon=0.5$ and $\epsilon \downarrow 0$) which is piecewise step-function, see Corollary~\ref{corr:x} in Appendix~\ref{app:new}. Along the red axis $\{(x,0,0)\;|\; x\in \R\}$ we have
$x=d(g,e)$. Bottom: the corresponding graph of $x \mapsto W_{n+1}^{\epsilon}((x,0,0),r_{n+1})$. As $n$ grows the staircase grows, as $\epsilon \to 0$ the size of the steps in the staircase vanishes and we see
$W^{\infty}(g)=d(g,e)$ in the most right column.
\label{fig:hot}}
\end{figure}

Details on the limit (\ref{stationary}), which takes place in the continuous setting before numeric discretization is applied, can be found in Appendix~\ref{app:new}. In Appendix~\ref{app:new} we provide solutions of (\ref{thePDE}) by a time shift in combination
with a morphological convolution\footnote{In fact, an `erosion' according to the terminology in morphological scale space theory, see e.g.~\!\cite{burgeth}.} with the corresponding morphological kernel, and show why the double limit is necessary. A quick intuitive explanation is given in Figure~\ref{fig:hot}, where we see that for $\epsilon>0$ we obtain stair-casing (due to a discrete rounding of the distance/value function) and where in the limit $\epsilon \downarrow 0$ the solution $W^{\infty}(g)=W(g)=d(g,e)$ is obtained.
\begin{remark}
The choice of our initial condition in Eq.~\!(\ref{thePDE1}) comes from the relation between linear and morphological scale spaces \cite{akian,burgeth}. Here, for linear $SE(2)$-convolutions over the $(\cdot,+)$-algebra
one has $\delta_{e}*_{SE(2)}U = U$.
For morphological $SE(2)$-convolutions (erosions) over the $(\min,+)$-algebra  \cite{DuitsJMIV} one has a similar property:
\begin{equation} \label{dilation}
(\delta^{M}_{e} \ominus U)(g):=\inf \limits_{q \in SE(2)}\left\{ \delta_{e}^{M}(q^{-1}g) + U(q)\right\} = U(g),
\end{equation}
This is important for representing viscosity solutions of left-invariant HJB-equations on $SE(2)$ by Lax-Oleinik \cite{Evans} type of formulas (akin to the $SE(3)$-case \cite{DuitsJMIV2013}).
\end{remark}

\begin{remark}
The stair-casing limit depicted in Figure~\ref{fig:hot} is similar to the basic Eikonal BVP on $\R$ with solution $d(x,0)=|x|$. On $\R$ the approach (\ref{thePDE}),(\ref{thePDE1})~and~(\ref{stationary}) provides pointwise limit {\small $|x|=\lim \limits_{\epsilon \to 0} \sum \limits_{m=0}^{\infty} r_{m+1} \, 1_{[r_m,r_{m}+\epsilon]}(|x|)=\lim \limits_{\epsilon \to 0} \sum \limits_{m=0}^{\lceil \frac{|x|}{\epsilon}\rceil} r_{m+1} \, 1_{[r_m,r_{m}+\epsilon]}(|x|)$}.
\end{remark}
\begin{remark}
By general semigroup theory \cite{akian}, one cannot impose both the initial condition and a boundary condition $W^{\epsilon}(e,r)=0$ at the same time, which forced us to update the initial condition (at $g=e$) in our iteration scheme (\ref{thePDE}).
The separate updating with value $0$ for $g=e$ in Eq.~(\ref{thePDE}) is crucial for the convergence in Eq.~\!(\ref{stationary}).
\end{remark}

\subsection{Construction of the Non-uniform Cost}\label{sec:cost}
The cost should have low values on locations with high curve saliency, and high values otherwise. Based on image $f$ we define the cost-function $\delta \leq \mathcal{C} \leq 1$ via
\begin{equation} \label{cost}
\mathcal{C}(x,y,\theta) = \frac{1}{1+ \lambda \left|\frac{(\mathcal{A}_{3}^{2}U_f)_+(x,y,\theta)}{\|(\mathcal{A}_{3}^2U_f)_+\|_{\infty}}\right|^p},
\end{equation}
with $\delta=\frac{1}{1+\lambda}$, and where $U_f:SE(2) \to \mathbb{R}$ is a lift of the image, with $\|\cdot\|_{\infty}$ the sup-norm, and
{\small
\[
(\mathcal{A}_{3}^2U_f)_+(g) = \left\{
\begin{array}{ll}
(-\sin \theta \partial_{x} + \cos \theta \partial_{y})^2 U_f(x,y,\theta) & \textrm{ if }(-\sin \theta \partial_{x} + \cos \theta \partial_{y})^2 U_f(x,y,\theta)\geq 0 \\
0 & \textrm{else }
\end{array}
\right.
\]
}
is a vessel detector where we use spatially isotropic Gaussian derivatives \cite{Franken2009}. The vessel-detector gives responses only if there are convex variations orthogonal to the elongated-structures of interest in $U_f(x,y,\theta)$.

The lifting is done using real-valued anisotropic wavelets $\psi$:
\begin{equation}
U_f(\mathbf{x},\theta) = \int_{{\mathbb{R}^{\rm{2}}}} \psi ( R_\theta ^{ - 1} (\mathbf{y} - \mathbf{x}))\, f(\mathbf{y})d\mathbf{y}.
\end{equation}
In this work we use the real part of so-called cake wavelets~\cite{DuitsIJCV2007} to do the lifting. These wavelets have the property that they allow stable reconstruction and do not tamper data evidence before processing takes place in the $SE(2)$ domain. Other type of 2D wavelets could be used as well. In related work by P\'echaud et al. \cite{Pechaud} the cost $\mathcal{C}$ was obtained via normalized cross correlation with a set of templates.

In Eq.~(\ref{cost}) two parameters, $\lambda$ and $p$, are introduced.
Parameter $\lambda$ can be used to increase the contrast in the cost function. E.g., by choosing $\lambda = 0$ one creates a uniform cost function and by choosing $\lambda> 0$ data-adaptivity is included. Parameter $p>1$ controls the steepness of the cost function, and in our experiments it is always set to $p=3$.
\begin{figure}[h]
\centering
\includegraphics[width=\textwidth]{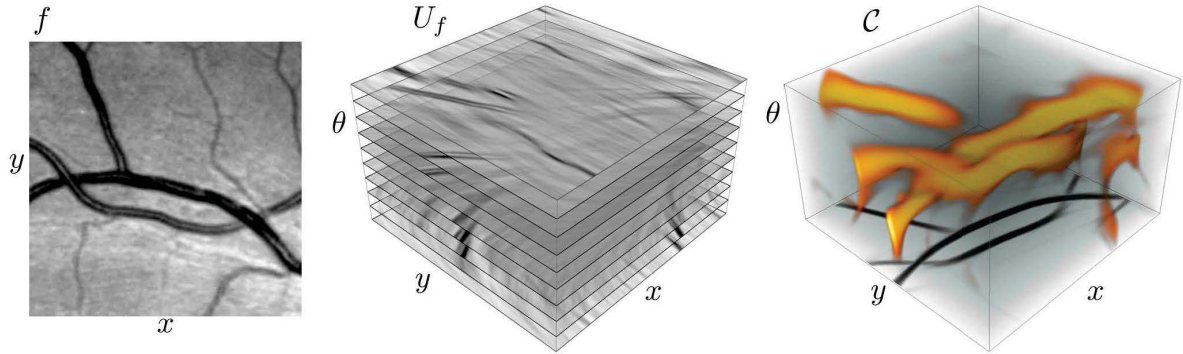}
\caption{Illustration of the cost function $\mathcal{C}$.
Left: retinal image patch $f$. Middle: corresponding function $U_f$ (`invertible orientation score') using the real part of a
cake-wavelet $\psi$ \cite{DuitsIJCV2007}. Right: the cost function $\mathcal{C}$ computed via Eq.~\!(\ref{cost}) visualized via volume rendering. The orange corresponds to locations where $\mathcal{C}$ has a low value. \label{figCost}}
\end{figure}

\section{Implementation} \label{ch:implementation}
To compute the SR geodesics with given boundary conditions we first construct the value function $W^{\infty}$ in Eq.~(\ref{Eik}), implementing the iterations at Eq.~(\ref{thePDE}), after which we obtain our geodesic $\gamma$ via a gradient descent on $W^{\infty}$ from $g$ back to $e$, recall Thm.~\ref{th:3a} (and Thm.~\ref{th:3}).
Throughout this section, we keep using the continuous notation $g \in SE(2)$ although within all numerical procedures $g$ is sampled on the following $(2 N+1) \times (2 N+1) \times (2 N_{\theta})$ equidistant grid:
{\small
\begin{equation} \label{grid}
\{(x_{i},y_j, \theta_k)\, | x_i = s_{xy} i, x_j= s_{xy} j,  \theta_{k}=k s_{\theta}, \textrm{ with } i,j=-N, \ldots N, k=-N_{\theta}\!+\!1,\ldots N_{\theta}\},
\end{equation}
}
with step-sizes  $s_{\theta}=\frac{2\pi}{2 N_{\theta}}$, $s_{xy}=\frac{x_{max}}{2N+1}$, with $N, N_{\theta} \in \mathbb{N}$. As a default we set $N=60$, $x_{max}=7$, $N_{\theta}=32$. The time-discretization grid is also chosen to be equidistant with time steps $\Delta r=\epsilon$.

On this grid we compute an iterative upwind scheme to obtain the viscosity solution $W^{\epsilon}$ at iteration Eq.~(\ref{thePDE}). Here we initialize $W^{\epsilon}(\cdot,0) = \delta_e^{MD}(\cdot)$, with the discrete morphological delta, given by $\delta^{MD}(g)=0 \textrm{ if }g=e$ and $1$ if $g\neq e$, and iterate
\begin{equation}\label{We}
\left\{
\begin{array}{lll}
    W^{\epsilon}(g,r + \Delta r) & = W^{\epsilon}(g,r) -  \Delta r \; H^{free}_D(g,dW^{\epsilon}(g,r)) & \mbox{for $g \neq e$}\\
    W^{\epsilon}(e,r + \Delta r) & = 0,                                                & \
\end{array}
\right.
\end{equation}
with free-time Hamiltonian (see Appendix~\ref{app:A}, Eq.~\!(\ref{the true})) given by
\[
H^{free}_D(g,dW^{\epsilon}(g,r)) = \left(   \frac{1}{\mathcal{C}(g)} \sqrt{   \beta^{-2} (\mathcal{A}_1 W^{\epsilon}(g,r))^2   +   (\mathcal{A}_2 W^{\epsilon}(g,r))^2   } - 1 \right),
\]
until convergence. We set $\Delta r=\epsilon$ in Eq.~\!(\ref{thePDE}).
In the numerical upwind scheme, the left-invariant derivatives are calculated via
\[
(\mathcal{A}_i W^{\epsilon}(g,r))^2 =
 \left(\max \left\{ \mathcal{A}_{i}^{-}W^{\epsilon}(g,r) , -\mathcal{A}_{i}^{+}W^{\epsilon}(g,r),0\right\}\right)^2,
\]
where $\mathcal{A}_i^+$ and $\mathcal{A}_i^-$ denote respectively the forward and backward finite difference approximations of $\mathcal{A}_i$.
Note that $W^{\epsilon}$ in (\ref{We}) is a first order finite difference approximation of $W^{\epsilon}_{n+1}$ in (\ref{thePDE}) at time interval
$r \in [n \epsilon, (n+1)\epsilon]$ and we iterate until the subsequent $\mathbb{L}_{\infty}$-norms differ less than $10^{-6}$.
This upwind scheme is a straightforward extension of
the scheme proposed in \cite{Rouy} for HJB-systems on $\R^{n}$. It produces sharp ridges at the 1st Maxwell set (cf.~\!Fig.~\!\ref{fig:gonz}) as it is consistent at local maxima.
For numerical accuracy and left-invariance we applied finite differences in the moving frame of left-invariant vector fields, using $B$-spline interpolation. This is favorable over finite differences in the fixed coordinate grid $\{x,y,\theta\}$. For details on these kind of left-invariant finite differences, and comparisons to other finite difference implementations (e.g. in fixed coordinate grid) see \cite{Franken2009}.

In our implementation the origin $e$ is treated separately as our initial condition is not differentiable.
We apply the update $W^{\epsilon}(e,r)=0$ for all $r \geq 0$.
We set step size $\epsilon = 0.1 \, \operatorname{min}(s_{xy} \beta, s_{\theta})$ with $s_{xy}$ and $s_{\theta}$ step sizes in respectively the $x$-$y$-directions and $\theta$-direction.

\section{Experiments and Results\label{sec:Experiments}}

\subsection{Verification for the Uniform Cost Case\label{ch:expcusp}}
Throughout the paper we have illustrated the theory with figures obtained via our new wavefront propagation technique. In this subsection we go through the figures that support the accuracy of our method. As the problem (\ref{problemP}) for $\mathcal{C}=1$ was solved \cite{yuriSE2FINAL,DuitsJMIV}, we use this as a golden standard for comparison.
\subsubsection{Comparison of BVP Solutions and the Cuspsurface}
\begin{figure}[t]
\centerline{
\includegraphics[width=\textwidth]{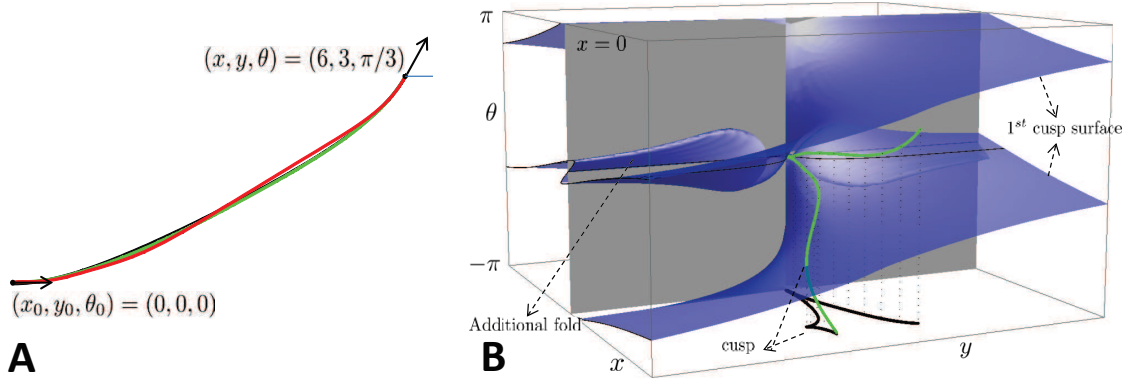}
}
\caption{
\textbf{A}: SR-geodesic example for the uniform cost case shows or PDE-discretizations (with 12 and 64 sampled orientations in red and green respectively) are accurate in comparison to analytic solutions in \cite{yuriSE2FINAL,DuitsJMIV} (in black).
\textbf{B}: The blue surface represents the cusp surface numerically computed via the proposed HJB-system (with $\mathcal{C}=1$) and subsequent calculation of the zero-crossings of $\mathcal{A}_1 W^{\infty}(x,y,\theta)$. Indeed if a SR-geodesic (in green) passes this surface, it passes in $\theta$-direction (with infinite curvature \cite{boscain,DuitsJMIV}), yielding a cusp on the spatial ground plane. The same blue surface is computed
in \cite[Fig.~11]{DuitsJMIV}. We even see the additional fold (top left passing the grey-plane) as some globally optimal SR-geodesics even exhibit 2 cusps.
\label{fig:cuspExample}
}
\end{figure}

Let us consider Fig.~\ref{fig:cuspExample}A. Here an arbitrary SR-geodesic between the $SE(2)$ points $\gamma(0)=e$ and $\gamma(T)=(6,3,\pi/3)$ is found via the IVP in \cite{yuriSE2FINAL} with end-time $T=7.11$ and initial momentum
\[
p_0= h_1(0) {\rm d}x + h_2(0) {\rm d}y  + h_3(0) {\rm d}\theta,
\]
with $h_{1}(0)=\sqrt{1-|h_{2}(0)|^2}$, $h_{2}(0)= 0.430$ and $h_{3}(0)= -0.428$, is used for reference (in black in Fig.~\ref{fig:cuspExample}A  .
Using the semi-analytic approach for solving the BVP in \cite{DuitsJMIV} an almost identical result is obtained.
The curves computed with our method with angular step-sizes of $2\pi/12$ and  $2\pi/64$ are shown in Fig.~\ref{fig:cuspExample}A in red and green respectively. Already at low resolution we observe accurate results. In Fig.~\ref{fig:gonz} we compare one SR-sphere  for $T=4$ (Fig.~\ref{fig:gonz}A) found via our method with the SR-wavefront departing from $e$ (Fig.~\ref{fig:gonz}B) computed by the method of characteristics \cite{yuriSE2}. We observe that our solution is non-smooth at the 1st Maxwell set $\mathcal{M}$ (\ref{Maxwell}) and that the unique viscosity solution stops precisely there, confirming Theorem~\ref{th:3}.
Finally, the blue surface in Fig.~\ref{fig:cuspExample}B represents the cusp surface, i.e. the surface consisting of all cusp points. Cusps are points that can occur on geodesics when they are projected into the image plane (see Fig.~\ref{fig:cuspExample}B). This happens at points $g_c$ where the geodesic is tangent to 
 $\left.\partial_\theta\right|_{g_c}=\left.\mathcal{A}_2\right|_{g_c}$. Then, the cusp surface $\gothic{S}_{cusp}$ is easily computed as a zero-crossing:
\begin{equation}\label{Scusp}
\begin{array}{l}
\gothic{S}_{cusp}:= \{g \in SE(2) \; |\; (\mathcal{A}_1 d)(g,e)=0\}, \
\gothic{S}_{cusp}^{num}:= \{g \in SE(2)\;|\; \mathcal{A}_1 W^{\infty}(g)=0\}.
\end{array}
\end{equation}
It is in agreement with the exact cusp surface $\gothic{S}_{cusp}$ analytically computed in \cite[Fig.~11]{DuitsJMIV}.
\begin{remark}
The simple key geometric idea behind (\ref{Scusp}) is that we have a cusp at time $t$ if
$u^{1}(t)=\frac{1}{\mathcal{C}^2(\gamma(t))} h_{1}(t)=\frac{1}{\mathcal{C}^2(\gamma(t))}\mathcal{A}_1 W^{\infty}(\gamma(t))=0$ which directly follows from Appendix~\ref{app:A}.
\end{remark}

\subsubsection{Comparison and Computation of SR-spheres}
Numerical verification of the solution obtained by our PDE approach was also done by comparison with the exact SR-distance, that was computed by explicit formulas for SR-geodesics (given on p.386 in~\cite{yuriSE2}) and the explicit formulas for cut time (that coincides with the first Maxwell time, given by (5.18)--(5.19) in~\cite{yuriSE2}). The experiments were done in the following way:
\begin{enumerate}
\item Compute a set of points:
$$EP(T) = \{(x_i,y_i,\th_i) = Exp(p_i, T)\;|\; p_i \in C, T \leq t_1^{MAX}(p_i), i =1,\ldots i_{max}\}$$
 of end points lying on the SR-sphere of fixed radius $T$ by analytic formulas for the exponential map and 1st Maxwell time $t_{1}^{MAX}$ ~\cite{yuriSE2}.
The number of end points $i_{max}$ was chosen as $i_{max} = 72\ T^2$.
Here $C$ is the cylinder in momentum space given by
\[
C=\left\{p \in T_{e}^{*}(SE(2))\;\right|\; \left. H^{fixed}(e,p)=1/2\right\},
\]
where we recall (\ref{Hfixed}).
The sampling points $p_i$ are taken by a uniform grid on the rectifying coordinates $(\varphi,k)$ of the mathematical pendulum (the ODE that arises in the PMP procedure, cf.\!~\cite[ch:3.2]{yuriSE2}), both for the rotating pendulum case ($p_{i} \in C_{2}$, yielding $S$-curves) and the oscillating pendulum case ($p_{i} \in C_{1}$, yielding $U$-curves), where we note that $C= \overline{C_{1} \cup C_{2}}$.
\item Evaluate the distance function $W^\infty(g_i)= W^\infty(x_i,y_i,\th_i)$ obtained by our numerical PDE-approach in section \ref{ch:implementation} for every point of the set $EP(T)$. We use 3rd order Hermite interpolation for $W^\infty(x_i,y_i,\th_i)$ at $g=g_i \in EP(T)$ in between the grid (\ref{grid}). 
\item Compute the maximum absolute error $\displaystyle{\max \limits_{g_i \in EP(T)}|W(g_i)-T|}$ and the maximum relative error $\displaystyle{\max \limits_{g_i \in EP(T)}(|W(g_i)-T|/T)}$.
\end{enumerate}
\begin{remark}
The exponential map $Exp: C \times \R^{+} \to SE(2)$ provides the end-point\\ $(x(t),y(t),\theta(t))=\gamma(t)=Exp(p_0,t)$ of the SR-geodesic $\gamma$, given SR-arclength $t$ and initial momentum $p_0 \in T_{e}^{*}(SE(2))$.
This exponential
map integrates the PMP ODE's in Appendix~\ref{app:A}.
\end{remark}

 The absolute and relative errors of the SR-distance computations for each end points located on SR-spheres of fixed radii are presented in Figure~\ref{fig:Tcompar}. The red graph corresponds to a sampling of $(N,N_{\theta})=(50,64)$, recall Eq.~(\ref{grid}), used in the SR-distance computation by our numerical PDE approach, and the blue graph corresponds to the finer sampling $(N,N_{\theta})=(140,128)$.  We see that the maximum absolute error does not grow, and that the relative error decreases when increasing the radius of the SR-sphere. Increase of sampling rate improves the result. For the finer sampling case, neither the absolute errors nor the relative errors exceed $0.1$.
\begin{figure}
\includegraphics[width=\hsize]{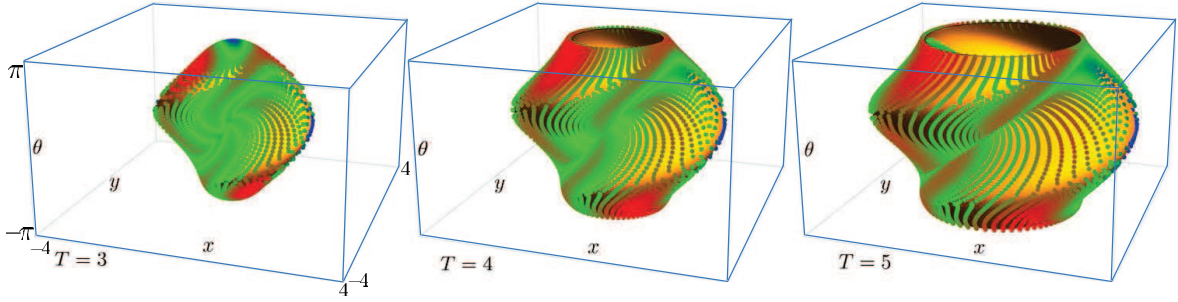}
\caption{Comparison of SR-spheres obtained by our numerical PDE-approach and the set of points $EP(T)$ lying on exact SR-spheres obtained by analytic formulas. From left to right: the SR-sphere with radius $t=T=3$, $T=4$ and $T=5$. The color indicates the difference between the exact and the numerical values of the SR-distance (blue for smallest, green for middle, and red for highest differences). Thus, we see  our algorithm is accurate, in particular along the fixed coordinate grid directions along $x$ and $\theta$-axis. }
\label{fig:colorpoints}
\end{figure}
\begin{figure}[t]
\centering
\begin{minipage}{0.35\textwidth}
\centering
\includegraphics[width=\textwidth]{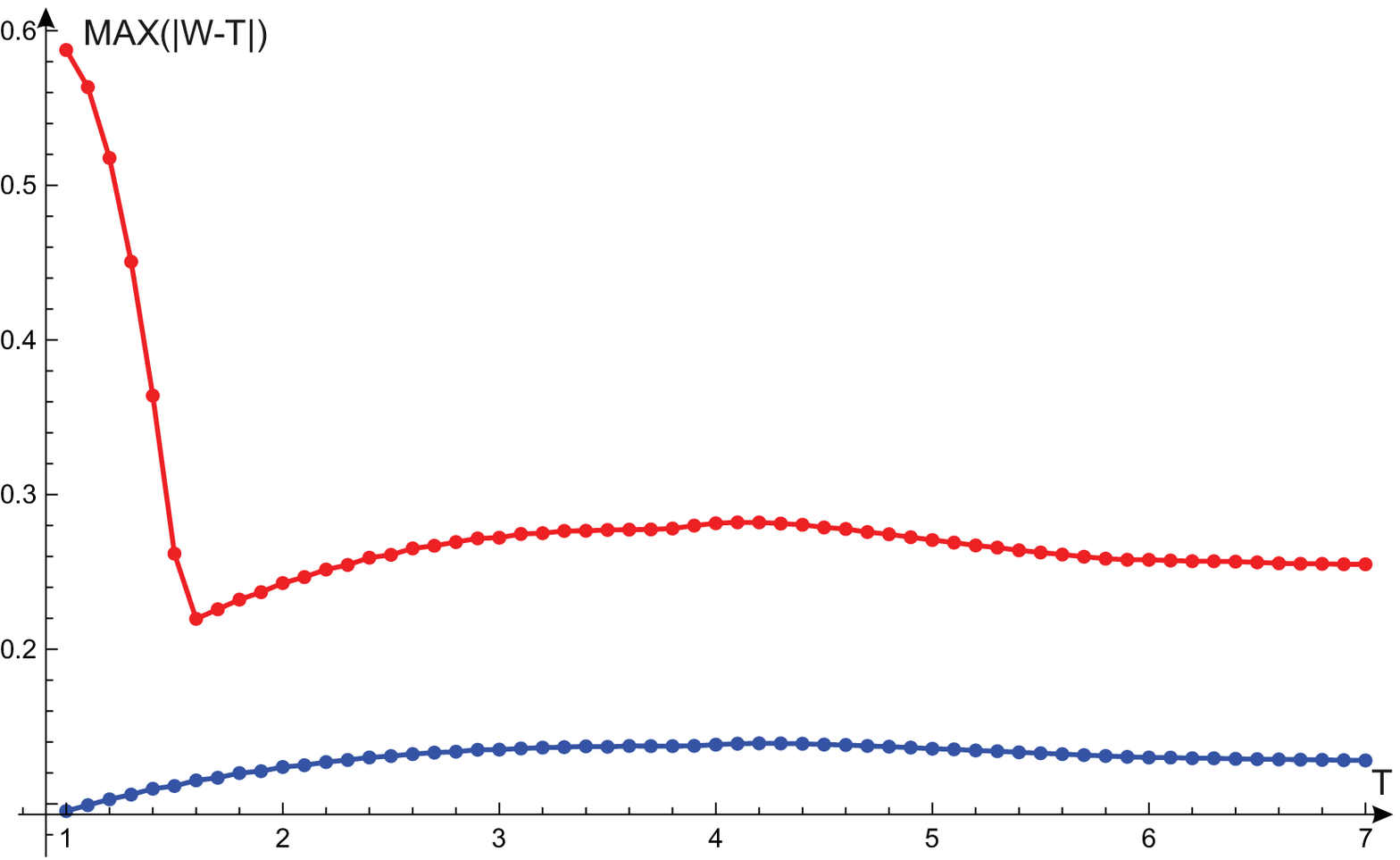}
A: Maximum absolute error
\end{minipage}
\begin{minipage}{0.35\textwidth}
\centering
\includegraphics[width=\textwidth]{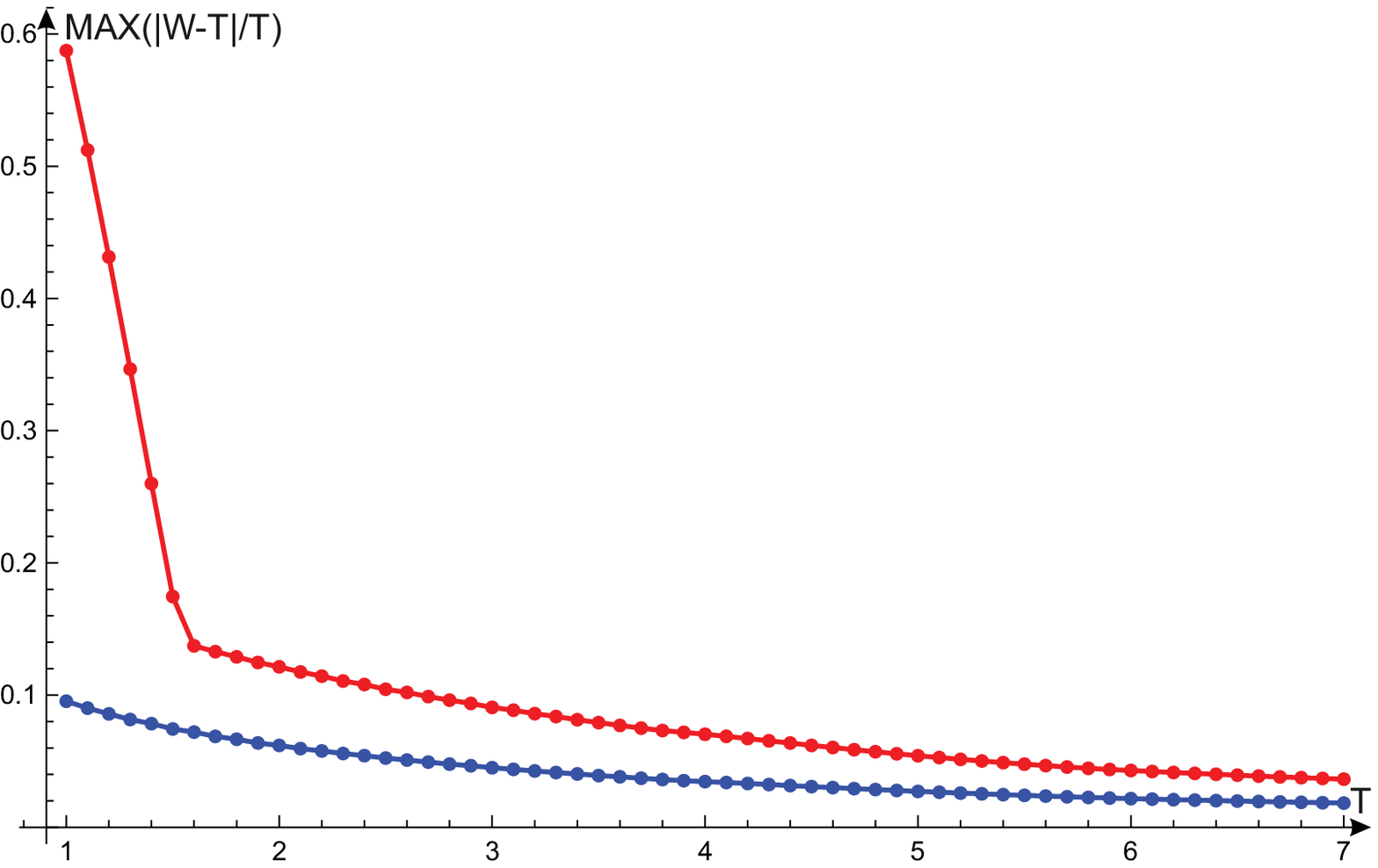}
B: Maximum relative error
\end{minipage}
\caption{Maximum error in computing of SR-distance for end points located on SR-spheres of different radii $t=T$ (from 1 to 7 with step 0.1), with number of end points $i_{max} = 72 \, T^2$. In red: errors are computed on a courser grid $(N,N_{\theta})=(50,64)$ and in blue: errors on a finer grid, with step sizes $s_{\theta}=\frac{2\pi}{N_{\theta}}$ and $s_{xy}=\frac{7}{N}$.}
\label{fig:Tcompar}
\end{figure}

\subsubsection{Comparison and Computation of 1st Maxwell Set}
\begin{figure}
\centerline{
\includegraphics[width=0.9\hsize]{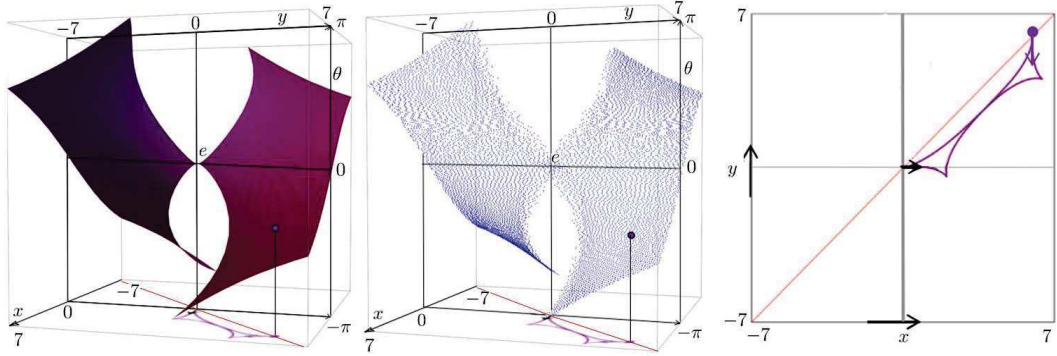}
}
\caption{Comparison of the 1st Maxwell set obtained by our numerical PDE-approach with the exact 1st Maxwell set~\cite{yuriSE2FINAL}.
Note that the local components of the 1st Maxwell set are part of the purple surface in Figure~\ref{globaloptim}.
\textbf{Left}: Local component of the exact Maxwell set $\mathcal{M}$ obtained by~\cite[thm 3.5]{yuriSE2FINAL} (where we recall that the cut locus coincides with the closure $\overline{\mathcal{M}}$ of the first Maxwell set~\cite[th:3.3]{yuriSE2CONJ}).
\textbf{Middle}: Local components of the Maxwell set $\mathcal{M}_{num}$ computed numerically by eq.~(\ref{eq:MaxwellFinDif}).
\textbf{Right}: Single case of a Maxwell point on the local part of the Maxwell set.
}
\label{fig:MaxwellSetsComparison}
\end{figure}

Considering forward and backward derivatives acting on the distance function $W^{\infty}(x,y,\theta)$, we can compute the 1st Maxwell set (recall eq.~(\ref{Maxwell}), see also Appendix~\ref{app:D}) as set of points where forward and backward left-invariant derivatives acting on distance function have different signs:
\begin{equation}
\label{eq:MaxwellFinDif}
\mathcal{M}_{num} = \displaystyle \bigcup \limits_{i=1}^{2} \{(x,y,\th) \in \SE |  \mathcal{A}_i^{+}W^{\infty}(x,y,\theta) > 0, \ \mathcal{A}_i^{-}W^{\infty}(x,y,\theta) < 0\}.
\end{equation}
Here $i=1$ corresponds to the \emph{local} component of the 1st Maxwell set (i.e. the part around $e$), and $i=2$ corresponds to the \emph{global} component of the 1st Maxwell set. The local component consists of two connected components lying on the surface given by $x\cos \frac{\theta}{2}+
y \sin \frac{\theta}{2}=0$ (i.e. the purple surface in Figure~\ref{globaloptim}), and the global component is a plane given by equation $\theta = \pi$ (for details, see~\cite{yuriSE2FINAL}).
In figure~\ref{fig:MaxwellSetsComparison} we compare the local component of $\mathcal{M}_{num}$ computed by our PDE approach with its exact counterpart $\mathcal{M}$, presented in~\cite[thm 3.5]{yuriSE2FINAL}. It shows that $\mathcal{M}_{num}$ is close to the exact $\mathcal{M}$, except for the conjugate points at the boundary (where (\ref{eq:MaxwellFinDif}) is less accurate).
Although note depicted here a similar picture was obtained for the global component, where $\mathcal{M}_{num}$ indeed covers the plane $\th = \pi$. Summaring, this experiment verifies the correctness of the proposed method, but it also shows that the method allows to observe the behavior of the 1st Maxwell set.
Eq.~\!(\ref{eq:MaxwellFinDif}) allows us to numerically compute the Maxwell set for the data-driven cases $\mathcal{C}\neq 1$ where exact solutions are out of reach.

\subsection{Feasibility Study for Application in Retinal Imaging\label{ch:appl}}

\begin{figure}
\centerline{
\includegraphics[width=0.85 \hsize]{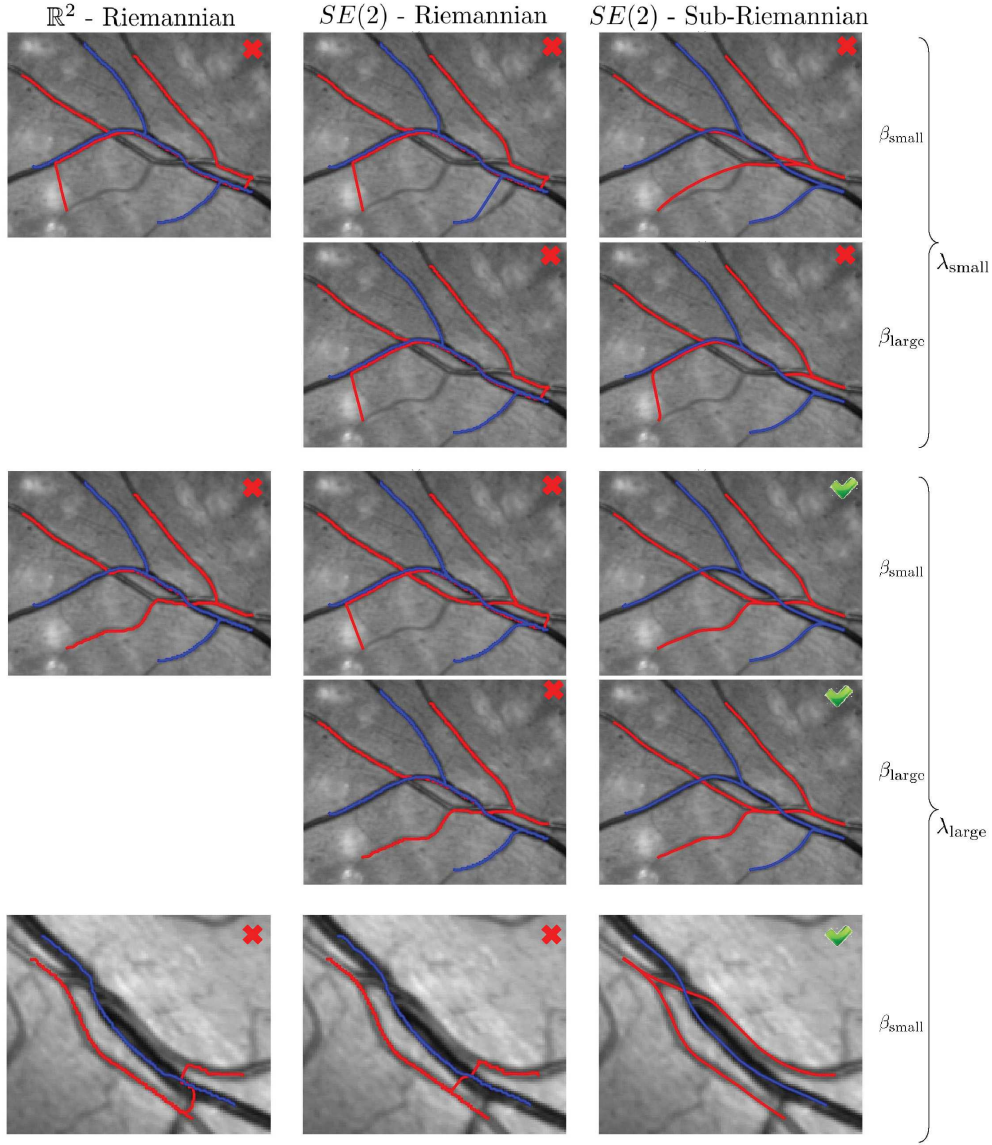}
}
\caption{
Data-adaptive sub-Riemannian geodesics obtained via our proposed tracking method  (Thm.~\ref{th:3}), with external cost (\ref{cost}),
with $p=3$, $\beta$ equals $\beta_{small}=0.01$, $\beta_{large}=0.1$ and $\lambda$ equals $\lambda_{small}=10$, $\lambda_{large}=100$.
We applied tracking from 2 seed-points each with several end-points (to test the crossings/bifurcations). To distinguish between tracks from the two seed-points we plotted tracts in different lighting-intensity. We indicated the valid cases only if all trajectories are correctly dealt with.
\label{fig:exp1}}
\end{figure}

\begin{figure}
\centerline{
\includegraphics[width=0.85 \hsize]{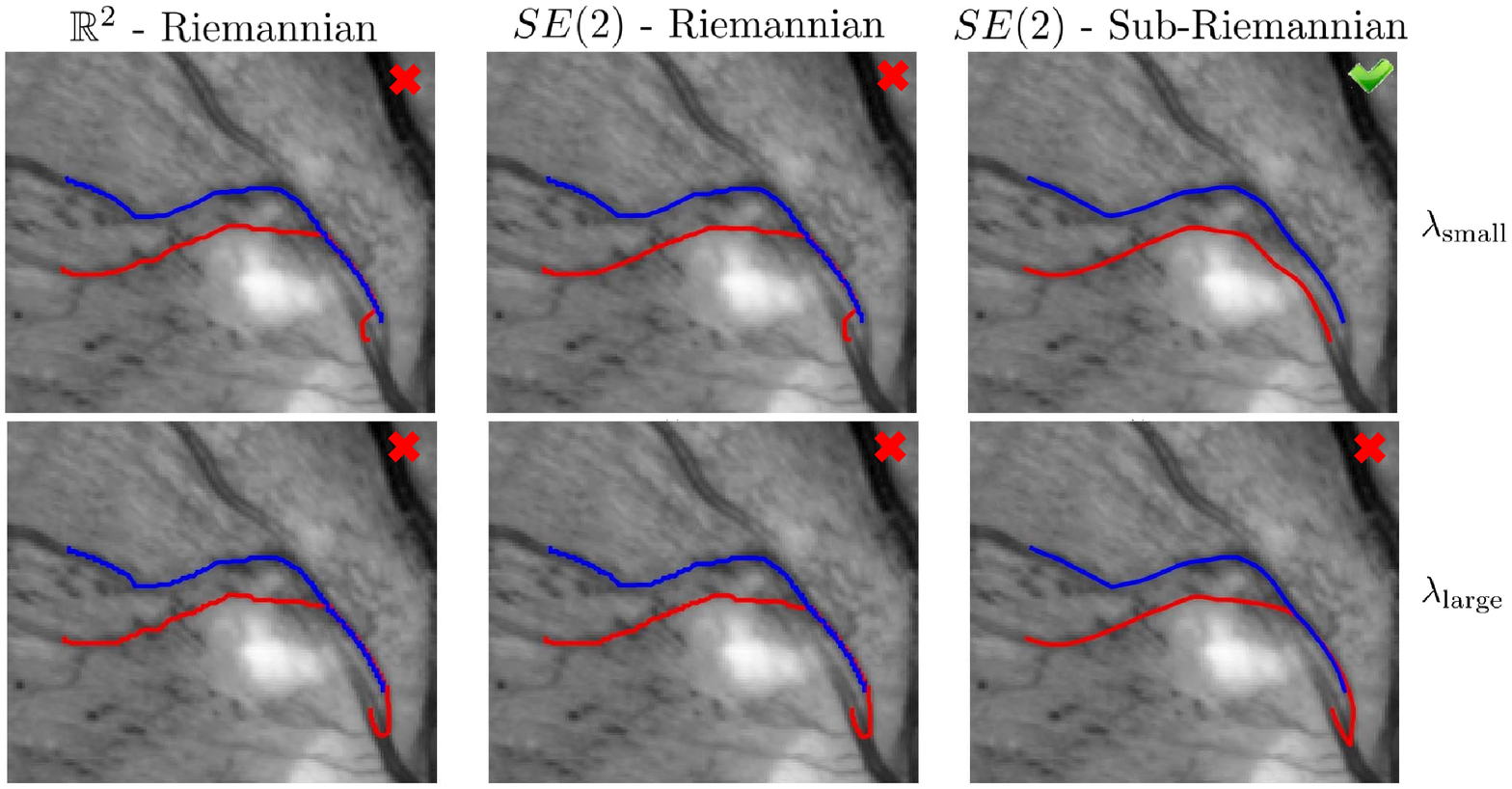}
}
\caption{Tractography (again for $\lambda=\lambda_{small}=10$, $\lambda=\lambda_{large}=100$ and $\beta=0.01$, $p=3$) in a patch of a challenging low-contrast retinal image of a diabetic patient.
In case of low-contrast (and less reliable cost) it is better to keep $\lambda$ small, in contrast to high contrast cases depicted in Figure~\ref{fig:exp1}. To distinguish between tracks from the two seed-points we plotted tracts in different lighting-intensity. We indicated the valid cases only if all trajectories are correctly dealt with.
\label{fig:badcase}}
\end{figure}

As a feasibility study for the application of our method in retinal images we tested the method on numerous image patches exhibiting both crossings, bifurcations, and low contrast, (Fig.~\!\ref{fig:exp1}, Fig.~\!\ref{fig:badcase}). 
For each seed point $g_0$ the value function $g \mapsto W^{\infty}(g_{0}^{-1}g)$ was calculated according to the implementation details in Section \ref{ch:implementation}, after which multiple end-points were traced back to the seed point.
The image dimensions of the patches is $180 \times 140$.

For the construction of the cost function (see e.g. Fig.~\ref{figCost}) we set $p=3$, and the lifting was done using cake wavelets with angular resolution $\pi/16$. More precisely we used a cake-wavelet with standard parameters $(N=8,N_{\theta}=32, s_{\theta}=\frac{\pi}{8}, \sigma_{s}=20px,\gamma=0.8)$, for details see cf.~\!\cite[ch:2]{BekkersJMIV}. The precise choice of anisotropic wavelet is not decisive for the algorithm (so other type of anisotropic wavelets and cost constructions could have been applied as well).

In all experiments we run with 4 settings for the two parameters $(\beta,\lambda)$ determining the sub-Riemannian geodesics, we
set $\beta_{small}=0.05$, $\beta_{large}=0.1$, $\lambda_{small}=10$, $\lambda_{large}=100$. The idea of these settings is to see the effect of the parameters, where we recall $\beta$ controls global stiffness of the curves, and $\lambda$ controls the influence of the external cost.
We also include comparisons to a Riemannian wavefront propagation method on $\R^{2}$, and a Riemannian wavefront propagation method on $SE(2)$.
These comparisons clearly show the advantage of including the sub-Riemannian geometry in the problem.
For results on two representative patches, see Figure~\ref{fig:exp1}. For results on 25 other patches see
\url{www.bmia.bmt.tue.nl/people/RDuits/Bekkersexp.zip}. Here, for fair and basic comparison of the geometries, we rely on the same cost function $\mathcal{C}$. That is,
we compare to
\begin{itemize}
\item Riemannian geodesics $\gamma(t)=(x(t),y(t),\theta(t))$ in $(SE(2), G_{full}^{\mathcal{C}})$ with
\[
\left.G_{full}^{\mathcal{C}}\right|_{\gamma(t)}(\dot{\gamma}(t),\dot{\gamma}(t))=(\mathcal{C}(\gamma(t)))^2\;
(|\dot{\theta}(t)|^2 + \beta^2 |\dot{x}(t)|^2 + \beta^2 |\dot{y}(t)|^2)
\]
\item
Riemannian geodesics $\ul{x}(s)=(x(s),y(s))$ in $(\R^{2}, G_{\R^2}^{\mathcal{C}})$ with metric tensor
\[
\left.G_{\R^{2}}^{\mathcal{C}}\right|_{\ul{x}(s)}(\dot{\ul{x}}(s),\dot{\ul{x}}(s)) = (c(\ul{x}(s)))^2 \, (|\dot{x}(s)|^2 +|\dot{y}(s)|^2),
\]
with $c(x(s),y(s))= \min \limits_{\theta \in [0,2\pi)} \mathcal{C}(x(s),y(s),\theta)$.
\end{itemize}
Typically, the wavefront propagation tracking methods on $(\R^{2},G_{\R^2})$ produces incorrect short-cuts at crossings and very non-smooth curves.
The Riemannian wavefront propagation tracking method (with spatial isotropy) $(SE(2),G_{full}^{\mathcal{C}})$ often deals correctly with crossings, but typically suffers from incorrect jumps towards nearly parallel neighboring vessels. Also it yields non-smooth curves.
This can be corrected for when including extreme anisotropy, see Remark~\ref{rem:RiemOK} below.
The Sub-Riemannian wavefront propagation method produces smooth curves that appropriately deals with crossings. For high contrast images with reliable cost $\mathcal{C}$ best results are obtained
with low $\beta$ and large $\lambda$. However, in low contrast images and/or patient data with severe abnormalities, low $\lambda$ is preferable, as in these cases the cost function is less reliable see Figure~\ref{fig:badcase}.
\begin{remark}\label{rem:RiemOK}
It is possible to construct a family of anisotropic Riemannian metric tensors, recall (\ref{dualbasis})):
$G^{\mathcal{C}}_{\epsilon}= \mathcal{C}^2 \, (\beta^{2}\omega^{1} \otimes {\omega}^{1}+ \omega^{2} \otimes \omega^{2}+ \beta^{2}\epsilon^2 \omega^{3}\otimes \omega^{3})$, which bridges the SR-metric $G^{\mathcal{C}}$ of our method (obtained by $\epsilon \to \infty$) to the full Riemannian metric tensor $G^{\mathcal{C}}_{full}$ (obtained by $\epsilon \to 1$).
For the values of $\beta$'s considered here, Riemannian geodesics and smooth Riemannian spheres
for highly anisotropic cases $\epsilon \geq 10$ approximate SR-geodesics and non-smooth SR-spheres.
In fact, with such extreme anisotropy in the Riemannian setting, the non-smooth ridges $\mathcal{M}$ in the SR spheres (see. e.g. the 1st Maxwell sets in Figure~\ref{fig:gonz}) are only little smoothed, and also the cusp-surface hardly changes. This observation allows to use the anisotropic fast-marching \cite{mirabeau} as an alternative fast method for computing the solution of (\ref{Eik}), instead of the iterative upwind finite difference approach in Section~\ref{sec:iterativeproc}.
\end{remark}

The experiments indicate $\beta=0.01$ (small) in combination with $\lambda=100$ (large) are preferable on our patches.
This typically holds for good quality retinal images of healthy volunteers. In lesser quality retinal images of diabetic patients, however, the cost function is less reliable and here $\lambda=10$ (small) can be preferable, see Figure~\ref{fig:badcase}.

However, it might not be optimal to set the $\beta$ parameter globally, as we did in these experiments, as smaller vessels are often more tortuous and therefore require more flexibility, see e.g.~\cite[fig.7]{ssvm}.
Furthermore, we do not include precise centerline extraction, which could e.g. be achieved by considering the vessel width as an extra feature (as in \cite{benmansour,Pechaud,LiYezzi}).

In future work we will pursue on sub-Riemannian fast-marching methods for fully automated vascular tree extraction starting from an automatically detected optic nerve head via state-of-the art method \cite{bekkersCMMVPR} followed by piecewise sub-Riemannian geodesics (comparable to  \cite{chen} in a Riemannian setting)
in between boundary points detected by a $SE(2)$-morphological approach. 
First experiments show that such a SR fast-marching method leads to considerable decrease computation time, hardly reduces the accuracy of the method,
and can be used to perform accurate and fast automatic full retinal tree segmentation. The advantage of such an approach over our previous work on automated vascular tree detection \cite{BekkersJMIV} is that each curve is a global minimizer of a formal geometric control curve optimization problem.
However, the SR-fast marching and automatic detection of the complete vascular tree via piecewise sub-Riemannian geodesics is beyond the scope of this theoretically oriented paper.


\section{Conclusion}

In this paper we propose a novel, flexible and accurate numerical method for computing solutions to the optimal control problem (\ref{problemP}), i.e. finding SR-geodesics in $SE(2)$ with non-uniform cost. The method consists of a wavefront propagation of geodesically equidistant surfaces computed via the viscosity solution of a HJB-system in $(SE(2),\Delta,G^{\mathcal{C}})$, and subsequent backwards integration, which gives the optimal tracks. We used PMP to derive both the HJB-equation and the backtracking. We have shown global optimality for the uniform case ($\mathcal{C}\!=\!1$) and that our method generates geodesically equidistant surfaces. Compared to previous works regarding SR-geodesics in $(SE(2),\Delta,G^{1})$ \cite{yuriSE2FINAL,DuitsJMIV,mashtakov}, we solve the boundary value problem without shooting techniques, using a computational method that always provides the optimal solution. Compared with wavefront propagation methods on the extended domain of positions and orientations in image analysis \cite{Pechaud,PechaudMICCAI}, we consider a SR metric instead of a Riemannian metric. Results in retinal vessel tracking are promising.

Fast, efficient implementation using ordered upwind schemes (such as the anisotropic Fast Marching method presented in \cite{mirabeau}) is planned as future work as well as adaptation to other Lie groups such as $SE(3)$, $SO(3)$. Of particular interest in neuroimaging is application to high angular resolution diffusion imaging (HARDI) by considering the extension to $SE(3)$ \cite{DuitsJMIV2013,PechaudMICCAI}.

\subsection*{Acknowledgements}

The research leading to these results has received funding from the European
Research Council under the European Community's Seventh Framework Programme
(FP7/2007-2013) / ERC grant \emph{Lie Analysis}, agr.~nr.~335555.
The authors gratefully acknowledge EU-Marie Curie project \emph{MANET}
ag.~no.607643 for financial support.

\appendix
\section{Application of PMP for Canonical Equations for Cost-adaptive Sub-Riemannian Geodesics}\label{app:A}
We study optimal control problem (\ref{problemP}). Recall Remark~\ref{rem:2}. Next we apply PMP to the action functional $J$ Eq.~\!(\ref{action}) with fixed total time $T>0$.
Since
$[\mathcal{A}_{i},\mathcal{A}_{j}]=\sum_{k=1}^{3}c^{k}_{ij}\mathcal{A}_{k}$, with non-zero coefficients $c^{3}_{12}=-c^{3}_{21}=-1$, $c_{23}^{1}=-c_{32}^{1}=-1$,
we have $[\Delta,\Delta]=T(SE(2))$ and we only need to consider normal trajectories.
Then the control dependent Hamiltonian of Pontryagin's Maximum Principle (PMP) expressed via left-invariant Hamiltonians $ h_i(p,g) = \langle p, \mathcal{A}_i(g) \rangle$, $i = 1,2,3$, with momentum $p \in T_g^{\ast}(SE(2))$, and $g = (x,y,\th) \in \R^2\times S^1$ reads as
\[H_u(p,g) = u^1 h_1(p,g) + u^2 h_2(p,g) - \frac12 \mathcal{C}^2(g)\left( \beta^2 |u^1|^2 +|u^2|^2\right).\]
Optimization over all controls produces the (maximized) Hamiltonian
\[H^{fixed}(g,p) = \frac{1}{2\, \mathcal{C}^2(g)}\left(\frac{h_1^2}{\beta^2} + h_2^2\right),\]
and gives expression for extremal controls
$u^1(t) = \frac{h_1(t)}{\mathcal{C}^2(\gamma(t)) \beta^2}, \quad u^2(t) = \frac{h_2(t)}{\mathcal{C}^2(\gamma(t))}$.
Using SR-arclength parametrization
$\left.G^{\mathcal{C}}\right|_{\gamma(t)}(\dot{\gamma}(t),\dot{\gamma}(t)) = 1$ implies $H^{fixed} = \frac12$ along extremal trajectories. We have the Poisson brackets
\begin{equation} \label{Poiss}
\begin{array}{l}
\{H,h_1\} = \frac{\mathcal{A}_1 \mathcal{C}}{\mathcal{C}} + \frac{h_2 h_3}{\mathcal{C}^2},\  \{H,h_2\} = \frac{\mathcal{A}_2 \mathcal{C}}{\mathcal{C}}  - \frac{h_1 h_3}{\beta^2 \mathcal{C}^2},\
\{H,h_3\} = \frac{\mathcal{A}_3 \mathcal{C}}{\mathcal{C}}  - \frac{h_2 h_1}{\mathcal{C}^2},
\end{array}\end{equation}
where $H=H^{fixed}$ and with
$\{F,G\}=\sum \limits_{i=1}^{3}
\frac{\partial F}{\partial h_i} \mathcal{A}_{i}G-
\frac{\partial G}{\partial h_i} \mathcal{A}_{i}F$. By Eq.~\!(\ref{Poiss}), by
$\{h_{i},h_j\}=
\mathcal{A}_{i}h_j-\mathcal{A}_{j}h_i=
\sum \limits_{k=1}^{3}c^{k}_{ij}h_k$, and by
$\dot{h}_i=\{H,h_i\}$, PMP gives us:
\begin{eqnarray}
p(\cdot) = \sum \limits_{i=1}^{3}h_{i}(\cdot) \left.\omega^{i}\right|_{\gamma(
\cdot)} \textrm{ and }
&&\begin{cases}
\dot{h}_1 = \frac{1}{\mathcal{C}(\gamma(\cdot))} \left.\mathcal{A}_1\right|_{\gamma(\cdot)} \mathcal{C} + \frac{h_2 h_3}{\mathcal{C}^2(\gamma(\cdot))},\\
\dot{h}_2 = \frac{1}{\mathcal{C}(\gamma(\cdot))} \left.\mathcal{A}_2\right|_{\gamma(\cdot)} \mathcal{C} - \frac{h_1 h_3}{\beta^2 \mathcal{C}^2(\gamma(\cdot))},\\
\dot{h}_3 = \frac{1}{\mathcal{C}(\gamma(\cdot))} \left.\mathcal{A}_3\right|_{\gamma(\cdot)} \mathcal{C} - \frac{h_2 h_1}{\mathcal{C}^2(\gamma(\cdot))},\\
\end{cases} \label{vertpart} \\
&&  \hspace{-2cm} \text{--- vertical part (for adjoint variables),} \nonumber \\
\dot{\gamma}(\cdot) = \sum \limits_{i=1}^{2}u^i(\cdot) \left.\mathcal{A}_{i}\right|_{\gamma(\cdot)} \textrm{ and }
&& \begin{cases}
\dot{x} = \frac{h_1}{\mathcal{C}^2(\gamma(\cdot)) \beta^2} \cos \th , \\
\dot{y} = \frac{h_1}{\mathcal{C}^2(\gamma(\cdot)) \beta^2} \sin \th , \\
\dot{\th} = \frac{h_2}{\mathcal{C}^2(\gamma(\cdot))},
\end{cases} \label{horpart} \\
&& \hspace{-2cm} \text{--- horizontal part (for state variables).} \nonumber
\end{eqnarray}
with dual basis $\{\omega^{i}\}$ for $T^{*}(SE(2))$ defined by $\langle \omega^{i}, \mathcal{A}_{j}\rangle =\delta_{j}^{i}$.

For a consistency check, we also apply the PMP-technique directly to Problem (\ref{problemP}) with free terminal time $T$, where
typically \cite{notes} the Hamiltonian vanishes. Then, using SR arclength parameter $t$,
the control dependent Hamiltonian of PMP equals
\[H_u(g,p) = u^1 h_1(p,g) + u^2 h_2(p,g) - \mathcal{C}(g)\sqrt{\beta^{2}|u^{1}|^2+|u^{2}|^2}.\]
Optimization over all controls
under SR arclength parametrization constraint \\
$\mathcal{C} \sqrt{\beta^2 |u^{1}|^2+|u^{2}|^2}=1$
produces via EL-optimization w.r.t. $(u^{1},u^2)$ (via unit Lagrange multiplier) the (maximized) Hamiltonian:
\begin{eqnarray}
\label{the true}
\begin{array}{ll}
H^{free}(g,p)&= \frac{1}{\mathcal{C}(g)}\sqrt{\frac{|h_{1}|^2}{\beta^2}+|h_{2}|^2}-1=0 \textrm{ with }p=\sum \limits_{i=1}^3h_{i}\omega^{i},
\end{array}
\end{eqnarray}
and by straightforward computations one can verify that both the horizontal part and the vertical part of PMP (but now applied to
$H^{free}$) is exactly the same as \eq{horpart} and \eq{vertpart}.
\begin{remark} The two approaches produce the same curves and equations, but their Hamiltonians are different. Nevertheless, we have
$H^{free}=0 \Leftrightarrow H^{fixed}=\frac{1}{2}$.
\end{remark}

\section{Lemmas Applied in the Proof of Theorem~\ref{th:3a}} \label{app:C}

In this section we consider preliminaries and lemmas needed for Thm.~\ref{th:3a}.
Before we can make statements on SR-spheres we need to explain the notion of geodesically equidistant surfaces, and their connection to HJB-equations. In fact, propagation of geodesically equidistant surfaces in
$(SE(2),\Delta,G^{\mathcal{C}})$ is described by a HJB-system on this SR-manifold.

Recall Remark~\ref{rem:2}.
Also recall, that in Appendix~\ref{app:A} we have applied PMP to this problem
yielding constant Hamiltonian $H^{fixed}=\frac{1}{2\mathcal{C}^{2}} (\beta^{-2} h_{1}^2 +h_{2}^2)=\frac{1}{2}$ relating to $H^{free}=\frac{1}{\mathcal{C}}\sqrt{\beta^{-2}h_{1}^2+h_{2}^{2}}-1=0$ via
$H^{free}= \sqrt{2 H^{fixed}}-1$.

In our analysis of geodesically equidistant surface propagation we first resort to the non-homogenous viewpoint on the Lagrangian and Hamiltonian (with fixed time), and then obtain the results on the actual homogeneous problem
(with free time) via a limiting procedure.
\begin{definition}\label{def:geodS}
Given $V:SE(2) \times \R^{+} \to \R$ continuous.
Given a Lagrangian  $L(\gamma(r),\dot{\gamma}(r))$ on the SR manifold $(SE(2),\Delta,G^{\mathcal{C}})$,
with $L(\gamma,\cdot): \Delta \to \R^{+}$ convex.
Then
the family of surfaces
\begin{equation} \label{geodsurface}
\mathcal{S}_{r}:=\{ g \in SE(2)
\;|\; V(g,r)=W_0(r)\}, \textrm{ with }
\end{equation}
$W_0:\R \to \R$ monotonic, smooth,
is geodesically equidistant if \mbox{
$L(\gamma(r),\dot{\gamma}(r))= W_{0}'(r)$.}
\end{definition}
\begin{remark} The motivation for this definition is
\[
\frac{d}{dR} \int \limits_{0}^{R} L(\gamma(r),\dot{\gamma}(r))\, {\rm d}r=
L(\gamma(R),\dot{\gamma}(R))=\frac{dW_0}{dr}(R).
\]
\end{remark}
\begin{lemma}\label{th:1}
Let $L$ be non-homogeneous and
$\lim \limits_{|v| \to \infty}\frac{L(\cdot,v)}{|v|} =\infty$.
Then the family of surfaces
$\{\mathcal{S}_{r}\}_{r \in \R}$ is geodesically equidistant if and only if
$V$ satisfies the HJB-equation
(where $r$ may be monotonically re-parameterized):
\begin{equation} \label{VPDE}
\begin{array}{l}
\frac{\partial V}{\partial r}(g,r)=
- H({\rm d}^{SR}V(g,r)), \textrm{ with }
{\rm d}^{SR}V(g,r)=\mathbb{P}_{\Delta}^{*} {\rm d}V(g,r)= \sum \limits_{i=1}^{2} \mathcal{A}_{i}V(g,r)\, \left.\omega^{i}\right|_{g}.
\end{array}
\end{equation}
Here $\mathbb{P}_{\Delta}^*(p=\sum \limits_{i=1}^{3} h_i\, \omega^{i})=
\sum \limits_{i=1}^{2} h_i\, \omega^{i}$ is a dual projection expressed in dual basis $\omega^{i}$ given by $\langle\omega^{i},\mathcal{A}_{j}\rangle =\delta^{i}_j$, and Hamiltonian
$H(\gamma,p)= \max \limits_{v \in T_{\gamma}(SE(2))} \{\langle p, v\rangle - L(\gamma,v)\}$.
\end{lemma}
\ \\
\textbf{Proof }
Substitute an arbitrary transversal minimizer $\gamma(r)$ into $V(\cdot,r)$ and
take the total derivative w.r.t. $r$:
\[
\frac{d}{dr} V(\gamma(r),r)=
\frac{\partial}{\partial r} V(\gamma(r),r)
+
\langle \left.{\rm d}V\right|_{\gamma(r)}, \dot{\gamma}(r) \rangle.
\]
Now $\gamma(r)$ on $S_{r}$, with tangent
$\dot{\gamma}(r)=\sum_{i=1}^{2}u^{i}(r) \left.\mathcal{A}_{i}\right|_{\gamma(r)}$, and thereby
\[
\frac{d}{dr} V(\gamma(r),r)=L(\gamma(r),\dot{\gamma}(r))=
\frac{\partial}{\partial r} V(\gamma(r),r)
+ \sum \limits_{i=1}^{2} u^{i}(r) \left.\mathcal{A}_{i}\right|_{\gamma(r)}V(\gamma(r),r).
\]
As a result we have
{\small
\begin{equation}\label{four}
\begin{array}{c}
-L(\gamma(r),\dot{\gamma}(r))+
\sum \limits_{i=1}^{2} u^{i}(r) \left.\mathcal{A}_{i}\right|_{\gamma(r)}V(\gamma(r),r) =
-\frac{\partial V}{\partial r} (\gamma(r),r)\overset{(1)}{\Leftrightarrow} \\
\sup \limits_{(u^{1}(r),u^2(r))\in \R^{2}}
\sum \limits_{i=1}^{2} u^{i}(r) h_i(r)-L(\gamma(r),\dot{\gamma}(r))
=
-\frac{\partial V}{\partial r} (\gamma(r),r) \overset{(2)}{\Leftrightarrow} \\
 H(\mathbb{P}_{\Delta}^{*}{\rm d}V(\gamma(r),r))=-\frac{\partial V}{\partial r} (\gamma(r),r),
\end{array}
\end{equation}
}
with components $h_i(r)=\left.\mathcal{A}_{i}\right|_{\gamma(r)}V(\gamma(r),r)$ of projected momentum vector
\[
\mathbb{P}_{\Delta}^{*}\, p(r)=\sum \limits_{i=1}^{2}h_i(r) \left.\omega^{i}\right|_{\gamma(r)} = \mathbb{P}_{\Delta}^{*}{\rm d}V(\gamma(r),r).
\]
Now every point $g \in S_{r}$ is part of a transversal minimizing curve $\gamma(r)$ and the result follows.
So the ``$\Rightarrow$'' is proven. Conversely, if the HJB-equation is satisfied it follows by the same computations (in reverse order)
that $L(\gamma(r),\dot{\gamma}(r))=\frac{d}{dr} V(\gamma(r),r)$, which equals $W_{0}'(r)$. 
$\hfill \Box$
\begin{remark}
In PMP \cite{notes} (see also Appendix~\ref{app:A}) the controls are optimized to obtain the Hamiltonian $H$ from the control dependent
Hamiltonian $H_u$. The first equivalence in (\ref{four}) is due to the maximum condition of PMP. The second equivalence in (\ref{four}), is by definition
of the Hamiltonian, where by the convexity assumption of the Lagrangian the supremum is actually a maximum \cite[ch:8]{Evans}.
\end{remark}

Next we apply the limiting procedure to obtain HJB-equations for geodesically equidistant surfaces in the actual homogeneous case of interest.
The actual \emph{homogeneous} Lagrangian case with $T$-free can be obtained
as a limit ($1\leq \eta \to \infty$) from non-homogeneous Lagrangian cases:
\begin{equation} \label{Leta}
L_{\eta}(\gamma(t),\dot{\gamma}(t))=
\frac{2\eta-1}{2\eta} \left(\left.G^{\mathcal{C}}\right|_{\gamma(t)}(\dot{\gamma}(t),\dot{\gamma}(t))\right)^{\frac{\eta}{2\eta-1}},
\end{equation}
and corresponding Hamiltonian (see Remark~\ref{rem:FT} below) equals
\begin{equation}\label{Heta}
H_{\eta}(\gamma(t),p(t))=
\frac{1}{2\eta} \left( \beta^{-2}h_{1}^2+h_{2}^2\right)^{\eta}
|\mathcal{C}(\gamma(t))|^{-2 \eta},
\end{equation}
and setting $r=t=W_0(t)$. Thus
$\frac{\partial V}{\partial r}(\gamma(r),r)=\frac{\partial V}{\partial t}(\gamma(t),t)=W_{0}'(t)=
L(\gamma(t),\dot{\gamma}(t))=\sqrt{\left.G^{\mathcal{C}}\right|_{\gamma(t)}(\dot{\gamma}(t),\dot{\gamma}(t))}=1$ in Eq.~\!(\ref{four}). Next we replace $V$ by $W$ to distinguish between the homogeneous and the non-homogeneous case).
\begin{lemma}\label{th:2}
The family of surfaces given by Eq.\!~(\ref{geodsurface}) is geodesically equidistant w.r.t. homogeneous Lagrangian
$L_{\infty}(\gamma,\dot{\gamma})=
\sqrt{\left.G^{\mathcal{C}}\right|_{\gamma}
(\dot{\gamma},\dot{\gamma})}$, with $r=t=W_0(t)$,
iff $W$ satisfies
 the HJB-equation
 \begin{equation} \label{Eik2}
\frac{1}{\mathcal{C}}\sqrt{\beta^{-2}|\mathcal{A}_{1}W|^2 +|\mathcal{A}_{2}W|^2}= 1 \desda H=0
 \end{equation}
where $H=\lim \limits_{\eta \to \infty} H_{\eta}=H^{free}$ the vanishing free-time Hamiltonian in Appendix~\ref{app:A}. Defining
Hamiltonian $\tilde{H}$ by
\begin{equation}\label{tildeH}
\tilde{H}(\gamma,p):=\mathcal{C}^{-1}(\gamma)\sqrt{\beta^{-2}h_{1}^2+h_{2}^2}
\end{equation}
puts
Eq.~(\ref{Eik2}) in Eikonal form $\tilde{H}({\rm d}^{SR}W(g,t))= 1$.
\end{lemma}
\ \\
\textbf{Proof }Tangential to the proof of Lemma~\ref{th:1}.
For $1\leq\eta< \infty$ we can apply Lemma~\ref{th:1} to Lagrangian $L_{\eta}$ given by (\ref{Leta}) whose associated Hamiltion
$H_{\eta}$ is given by (\ref{Heta}) due to PMP (or just the Fenchel transform on $\R^{2}$). In the limiting case $\eta \to \infty$, where the Lagrangian is homogeneous and the Hamiltonian vanishes. Finally we note that now we have
\[
\frac{\partial W}{\partial r}(\gamma(r),r)=
\frac{\partial W}{\partial t}(\gamma(t),t)= W_{0}'(t)= L(\gamma(t),\dot{\gamma}(t))= 1,
\]
from which the result follows.  $\hfill \Box$
\begin{remark}
The relation between the various Hamiltonians is
\[
H_{\eta \to \infty}=H^{free}=\sqrt{2 H^{fixed}}-1=\sqrt{2 \, H_{\eta=1}}-1=\tilde{H}-1=0.
\]
\end{remark}
\begin{remark} \label{rem:FT}
The relation between the Lagrangian $L_{\eta}$ given by (\ref{Leta}) and the Hamiltonian
(\ref{Heta}) is the (left-invariant, SR) Fenchel transform on $SE(2)$. Due to left-invariance this Fenchel-transform actually boils down to an ordinary Fenchel-transform on $\R^{2}$ when expressing velocity and momentum in the left-invariant frame. Indeed we have
\begin{multline}
H_{\eta}(\gamma,p) =
[\gothic{F}_{\mathcal{L}(SE(2)) \cap \Delta} (L_{\eta}(\gamma,\cdot))](p):=\\
  \sup \limits_{(u^{1},u^{2}) \in \R^{2}} \{ -\frac{2\eta-1}{2\eta}(\mathcal{C}(\gamma))^{\frac{2\eta}{2\eta-1}}
  (\beta^{2}|u^{1}|^2 +|u^{2}|^2)^{\frac{\eta}{2\eta-1}} + h_{1}u^{1}+h_{2}u^{2}\}
  \end{multline}
with horizontal velocity $v=u^{1}\mathcal{A}_{1} +u^{2} \mathcal{A}_{2}$,
and momentum $p=\sum_{i=1}^{3}h_{i}\omega^{i}$.
\end{remark}


\section{Viscosity Solutions for HJB-systems in SE(2)\label{app:B}}
\begin{definition} The (Cauchy problem) for a HJB-equation (akin to~\cite[ch:10.1]{Evans}) on $SE(2)$ is given by
\begin{equation} \label{oneH}
\begin{cases}\frac{\partial W}{\partial t} = - H(g,{\rm d}^{SR}W) \text{ in } SE(2) \times (0,T),
\\ W(g,0) = W_0, \end{cases}
\end{equation}
whereas a boundary value problem for HJB-equation is given as
\begin{equation} \label{twoH}
H(g,{\rm d}^{SR}W) = 0 \textrm{ on } SE(2) \setminus \{e\},\qquad  W(e)=0;
\end{equation}
where $T>0$ is prescribed, $W_0$ is a given function (or a cost measure \cite{akian}), $H(g,p) = H^{free}(g,p)$ is the free-time Hamiltonian given by (\ref{Hfree}), and ${\rm d}^{SR} W =  \sum \limits_{i=1}^{2} \mathcal{A}_{i}W(g,t) \, \left.\omega^{i}\right|_{g}$. 
\end{definition}

\begin{remark}
Combined Cauchy-Dirichlet problems exist \cite{trelat}, but they are defined on (analytic) open and bounded domains. Thereby they cannot be applied to our set of interest $SE(2)\setminus \{e\}$ as this would violate semigroup theory \cite{akian,Evans,yosida,burgeth}. This is also clear in view of the Cramer transform \cite{akian}, putting an isomorphism between HJB- and diffusion systems.
\end{remark}
\begin{remark}
In Eq.~\!(\ref{twoH}) it is crucial that the free time Hamiltonian is used. In the definition of viscosity solutions of the Cauchy problem, Eq.~\!(\ref{oneH}), one can set both $H=H^{free}$ (as done in the body of the article) or $H=H^{fixed}$ as done in Appendix~\ref{app:C}.
\end{remark}

HJB-systems in general do not have unique (weak) solutions.
To avoid multiple (non-desirable) solutions,
one must impose the viscosity condition \cite{Evans,Lions} commonly applied in wavefront methods acting directly in the image domain $\R^2$ \cite{osher,sethian}.
The viscosity solution is obtained by the vanishing viscosity method~\cite{Lions}.
The idea of this method is to add to the HJB-equation a term $\varepsilon \Delta$ and to pass to the limit, when $\varepsilon$ goes to 0. Here $\Delta$ denotes the Laplacian, that in our case (for $\mathcal{C}=1$) equals $\Delta^{SR} = \sum_{i=1}^{2} \mathcal{A}_{i}(\beta_{i})^{-2}\,\mathcal{A}_i$, with $\beta_1=\beta, \beta_2=1$.
Here the name ``viscosity solutions'' comes from fluid dynamics, where typically the term $\varepsilon \Delta$ represents a physical viscosity.
For an intuitive illustration of the geometric property of such solutions see \cite[fig.30]{bressan}.
 The viscosity solution of the initial value problem can be defined alternatively as follows. \\
 \
 \begin{definition}
 Let $H(g,\cdot)$ be a convex Hamiltonian for all $g \in SE(2)$ s.t. $H(g,p) \to \infty$ if $p\to \infty$. The function $W: SE(2) \times \R \to \R$ is viscosity solution of $\frac{\partial W}{\partial t} = - H(g,{\rm d}^{SR}W)$ if it is a weak solution such that for all functions $V \in C^{1}(SE(2)\times \R, \R)$ one has
\begin{itemize}
\item if $W-V$ attains a local maximum at $(g_0,t_0)$ then
$\left.\left(\frac{\partial V}{\partial t}  + H(g,{ \rm d}^{SR}V)\right)
\right|_{(g_0,t_0)}  \leq 0$,
\item if $W-V$ attains a local minimum at $(g_0,t_0)$ then
$\left.\left(\frac{\partial V}{\partial t}  + H(g,{\rm d}^{SR}V)\right)
\right|_{(g_0,t_0)}\geq 0$.
\end{itemize}
\end{definition}
Similarly, the viscosity solution of the boundary value problem (that is equivalent to Eikonal equation, when $t$ is SR-arclength) can be defined as follows:
\begin{definition} \label{def:VSEik}
A solution $W:SE(2) \to \R$ of Eq.~\!(\ref{twoH}) is called a viscosity solution if
for all functions $V \in C^{1}(SE(2), \R)$ one has
\begin{itemize}
\item if $W-V$ attains a local maximum at $g_0$ then
$H^{free}(g_0,{\rm d}^{SR}V)) \leq 0$,
\item if $W-V$ attains a local minimum at $g_0$ then
$H^{free}(g_0,{\rm d}^{SR}V)) \geq 0$.
\end{itemize}
\end{definition}

\section{Proof of Theorem \ref{th:3}~\label{app:D}}
\ \\
The back-tracking (\ref{steepest})
is a direct result of Lemma~\ref{th:2} in Appendix~\ref{app:C} and
PMP in Appendix~\ref{app:A}. According to these results
one must set
\[
u^1(t)=\frac{h_{1}(t)}{(C(\gamma(t)))^2 \beta^2}=
\frac{\left.\mathcal{A}_{1}W\right|_{\gamma(t)}}{(C(\gamma(t)))^2 \beta^2} \textrm{ and }
u^2(t)=\frac{h_{2}(t)}{(C(\gamma(t)))^2 }=
\frac{\left.\mathcal{A}_{2}W\right|_{\gamma(t)}}{(C(\gamma(t)))^2}
\]
from which the result follows. Then we recall from Thm.~\ref{th:3a} that
$\mathcal{S}_{t}$ given by (\ref{St}) are geodesically equidistant surfaces
propagating with unit speed from the origin. So $\mathcal{S}_{t}$ are
candidates for sub-Riemannian spheres, but it remains to be shown that the
back-tracking (\ref{steepest}) will neither pass a Maxwell point or
a conjugate point, i.e. $t\leq t_{cut}$. Here $t_{cut}$ denotes cut time, where a geodesic looses its optimality.

At Maxwell points $g^*$ induced
by the 8 reflectional symmetries \cite{yuriSE2} two distinct SR-geodesics meet with
the same SR distance. As SR-geodesics in $(SE(2),\Delta,G^{1})$ are analytic \cite{yuriSE2},
these two SR-geodesics do not coincide on an open set
around end-condition $g^{*}$, and the SR spheres are non-smooth at $g^{*}$.
Regarding the set $\mathcal{M}$, we note that the Maxwell sets related to the $i$-th reflectional symmetry $\epsilon_i$
are defined by
\[
\begin{array}{l}
\textrm{MAX}^{i}=\{(p_0,t) \in T_{e}^{*}(SE(2)) \times \R^{+} \;|\; H(p_0)=\frac{1}{2} \textrm{ and }Exp(p_0,t)=Exp(\epsilon_i p_0, t) \}, \\
\textrm{max}^{i}=Exp(\textrm{MAX}^{i}), i=1,\ldots 8,
\end{array}
\]
where we may discard indices $i=3,4,6$ as they are contained in $\{\max_{1},\max_{2},\max_{5},\max_7\}$.
Now with $\widetilde{\textrm{max}}^{i}$ we denote the Maxwell set with \emph{minimal} positive Maxwell times over all symmetries
(i.e. we include the constraint $t\leq \min\{t_{max}^{i}\}$ where the minimum is taken over all positive Maxwell times along each trajectory), then we find $\mathcal{M}$ to be contained within the union of the following sets\footnote{In \cite[Eq.3.13]{yuriSE2FINAL} it is shown that
$\widetilde{\textrm{max}}^{2}=\{(x,y,\theta) \in SE(2)\,| \,y \sin \theta/2+ x \cos \theta/2=0\textrm{ and }|-x \sin (\theta/2) + y \cos (\theta/2)|> |R_{1}^{1}(\theta)|\}$
with $R_{1}^{1}$ defined in \cite[Lemma 2.5]{yuriSE2FINAL}.
We also observed such a loss of the Maxwell point property in our numerical algorithm,
as kinks in $W(g)=t$ can disappear when moving on the set $y \sin \theta/2+ x \cos \theta/2=0$. See Figure~\ref{fig:MaxwellSetsComparison}.
}:
\[
\begin{array}{l}
\widetilde{\textrm{max}}^{2}\subset\{(x,y,\theta) \in SE(2)\;| \; y \sin \theta/2+ x \cos \theta/2=0\}, \
\widetilde{\textrm{max}}^{5}=\{(x,y,\theta) \in SE(2)\;| \;\theta= \pi \},
\end{array}
\]
where \cite[th:5.2]{yuriSE2} shows we must discard the first reflectional symmetry $\epsilon_{1}$ as it does not produce Maxwell points. Now
for generic geodesics (not passing the special conjugate points that are limit points of Maxwell points and not Maxwell points themselves) $t_{cut}= t_{MAX}^{1}$, as proven in \cite[th:3.3]{yuriSE2CONJ}, where $t_{MAX}^{1}>0$ denotes the first Maxwell time associated to the $8$ discrete reflectional symmetries.

\begin{figure}\centerline{
\includegraphics[width=0.62\hsize]{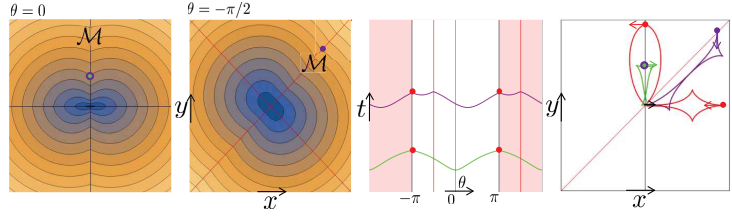}
}
\caption{Overview of Maxwell points. Two Maxwell points on the purple surface $x \cos \frac{\theta}{2} + y \sin \frac{\theta}{2}=0$ and two red Maxwell points on the surface $|\theta|=\pi$ (recall Figure~\ref{globaloptim}). In all cases we see
that local kinks in the viscosity solutions are upward, and the back-tracking algorithm can not pass these points.
 \label{fig:MaxwellDemo}}
\end{figure}

During the back-tracking the set $\mathcal{M}$ is never reached at internal times (only when starting at them, recall Remark~\ref{rem:MW}), since they are ``uphill'' from all possible directions during dual steepest descent tracking (\ref{steepest}), as we will shown next. As a result we have $t \leq t_{cut}=t_{MAX}^{1}$.
Consider Fig.~\ref{globaloptim}. At Maxwell points $g^* \in \mathcal{M}$ due to the reflectional symmetries there exist two distinct directions in the 2D-horizontal part $\Delta_{g^*}$ of the tangent space $T_{g^*}(SE(2))$ where the directional derivative is positive. If there would be a direction in the tangent space where the directional derivative is negative then
there would be a direction in $\Delta_{g^*}$ with zero directional derivative of $W(\cdot)$ at $g^{*}$ towards the interior of the sphere yielding contradiction. Here we note that due to the viscosity property of the HJB-solution, kinks at the Maxwell points are pointing upward (see Fig.~\ref{globaloptim} and Fig.~\ref{fig:MaxwellDemo}) in the backward minimization tracking process \cite[fig.30]{bressan}. Furthermore, we note that SR spheres $\mathcal{S}_{t}$ are continuous
\cite{yuriSE2FINAL} and compact, as they are the preimage $\mathcal{S}_{t} = d(\cdot,e)^{\leftarrow}(\{t\})$ of compact set $\{t\}$ under continuous mapping $d(\cdot,e)$.
Continuity of $d(\cdot,e)$ implies the spheres are equal to the 2D-boundaries of the SR balls (w.r.t. the normal product topology on $\R^{2}\times S^{1}$).

The algorithm also cannot pass conjugate points that are limits of 1st Maxwell points, but not Maxwell points themselves. See Fig.~\ref{fig:gonz}.
Such points exist on the surface $R_{2}=0$ and are by definition within $\overline{\mathcal{M}} \setminus \mathcal{M}$.
Suppose the algorithm would pass such a point at a time $t>0$
(e.g. there exist 4 such points on the sphere with radius $4$, see Fig.~\ref{globaloptim})
then due to the astroidal shape of the wavefront at such a
point, cf.~\!\cite[Fig.11]{yuriSE2CONJ}, there is a close neighboring tract that would pass a 1st Maxwell point which was already shown to be impossible (due to the upward kink-property of viscosity solutions). $\hfill \Box$

\begin{remark}
The sub-Riemmannian spheres are non-smooth only at the 1st Maxwell set $\mathcal{M}$. They are smooth at the conjugate points in
$\overline{\mathcal{M}}\setminus \mathcal{M}$ (where the reflectional symmetry no longer produce two curves/fronts). In the other points on $\mathcal{S}_{t} \setminus \overline{\mathcal{M}}$ the SR-spheres are locally smooth
(by the Cauchy-Kovalevskaya theorem and the semigroup property of the HJB-equations).
\end{remark}

\section{The Limiting Procedure (\ref{stationary}) for the sub-Riemannian Eikonal Equation\label{app:new}}

In this section we study the limit procedure (\ref{stationary}), illustrated in Figure~\ref{fig:hot}.
To this end we first provide a formal representation of the viscosity solutions of system (\ref{thePDE}),
where we rely on viscosity solutions of morphological scale spaces obtained by super-position over the $(\min,+)$ algebra, i.e. obtained by morphological
convolution (erosions) with the morphological impulse response, cf.~\!\cite{burgeth}. Now as the HJB-equations of such morphological scale spaces do not involve
a global offset by $1$ in the right-hand side of the PDE we need to combine such erosions with a time shift in order to take into account the global offset.
It turns out that the combination of these techniques provide staircases with steps of size $\epsilon$, so that we obtain the appropriate limit by taking
the limit $\epsilon \to 0$ afterwards as is done in (\ref{stationary}).

Morphological convolutions over the $SE(2)$ group are obtained by replacing in linear left-invariant convolutions (likewise the $SE(3)$-case studied in \cite{DuitsJMIV2013}), the usual $(+,\cdot)$-algebra by the $(\min,+)$-algebra.
Such erosions on $SE(2)$ are given by
\begin{equation} \label{erosion}
(k \ominus f)(g):= \inf \limits_{h \in SE(2)} \{k(h^{-1}g)+ f(h)\}.
\end{equation}
Furthermore, to include the updating of the initial condition in (\ref{thePDE}) we define
\begin{equation}\label{tilde}
\tilde{W}(g):=\left\{
\begin{array}{ll}
W(g) &\textrm{ if }g \neq e \\
0 &\textrm{if }g=e.
\end{array}
\right.
\end{equation}
\begin{lemma} \label{lemma:imp}
Let $\epsilon >0$, $n \in \mathbb{N}$.
The viscosity solution of (\ref{thePDE}) is given by
\begin{equation}\label{mainidea}
W_{n+1}^{\epsilon}(g,r) = (k_{r-n \epsilon} \ominus \tilde{W}_{n}^{\epsilon})(g) + (r -n \epsilon),
\end{equation}
for $r \in [r_n,r_{n+1}]=[n \epsilon, (n+1)\epsilon]$ and morphological kernel
$k_{v}(g)$, $v\geq 0$, given by
\[
k_{v}(g)=
\left\{
\begin{array}{ll}
0 &\textrm{if }d(g,e) \leq v, \\
\infty &\textrm{else,}
\end{array}
\right.
\]
where $d(g,e)$ denotes the Carnot-Caratheodory distance (\ref{SRdist}) between $g \in SE(2)$ and $e=(0,0,0)$.
For $n=0$ we have that the viscosity solution of (\ref{thePDE1}) is given by $W_{1}^{\epsilon}(g,r)=k_{r}(g)+r$.
\end{lemma}
\textbf{Proof }
In order to care of the constant off-set in the HJB-equatons of (\ref{thePDE}) and (\ref{thePDE1}), we set $r=r_{new}+r_n$ and we define for $n=0,1,2,\ldots$ the
functions $V^{\epsilon}_{n+1}:SE(2)\times [0,\epsilon] \to \R$ by
\[
V^{\epsilon}_{n+1}(g,r_{new}):= W^{\epsilon}_{n+1}(g,r_{new}+r_{n}) - r_{new},
\]
with $r_{new} \in [0,\epsilon]$ and $V^{\epsilon}_{n+1}$ the viscosity solution of
\[
\left\{
\begin{array}{l}
\frac{\partial V_{n+1}^{\epsilon}}{\partial r_{new}} (g, r_{new})= -1+1-\tilde{H}({ \rm d}^{SR} V_{n+1}^{\epsilon}(g,r_{new}))=-\tilde{H}({ \rm d}^{SR} V_{n+1}^{\epsilon}(g,r_{new})), \\
\textrm{for }g\neq e \textrm{ we have } V_{n+1}^{\epsilon}(g,0)=
\left\{
\begin{array}{ll}
\infty &\textrm{if }n=0, \\
W_{n}^{\epsilon}(g,r_n) &\textrm{if }n \in \mathbb{N}
\end{array}
\right.
\\
\textrm{for }g=e \textrm{ we have }V_{n+1}^{\epsilon}(g,0)=V_{n+1}^{\epsilon}(e,0)=0.
\end{array}
\right.
\]
where we use short notation for the sub-Riemannian derivative ${\rm d}^{SR}V:= \sum_{i=1}^{2}\mathcal{A}_{i}V \omega^{i}$, recall
(\ref{dualbasis}) in Remark~\ref{rem:3}, and where Hamiltonian $\tilde{H}$ is given by (\ref{tildeH}).

Now let us first consider the case $n=0$. By the results in Appendix~\ref{app:C} the Hamiltonian system (\ref{thePDE1}) provides geodesically equidistant wavefront propagation traveling with unit speed and departing directly from the unity element.
As a result, we find
\[V_{1}^{\epsilon}(g,r_{new})=k_{r_{new}}(g)=
\left\{
\begin{array}{cc}
0 & \textrm{ if }d(g,e)\leq r_{new} \\
\infty  & \textrm{ else. }
\end{array}
\right.
\]
and by left-invariant `superposition' over the $(\min,+)$-algebra we find for $n=1,2,\ldots$ that
$V_{n+1}^{\epsilon}(g,r_{new})= (k_{r_{new}} \ominus \tilde{W}_{n}^{\epsilon}(\cdot,r_{n}))(g)$. Finally, we have
\[
W_{n+1}^{\epsilon}(g,r)=V_{n+1}^{\epsilon}(g,r-n\epsilon)+ r-n \epsilon=
(k_{r-n \epsilon} \ominus \tilde{W}_{n}^{\epsilon}(\cdot,r_n))(g)+ r -n \epsilon.
\]
\begin{corollary} \label{corr:x}
Let $n \in \mathbb{N}$. Let $\epsilon>0$. The following identity holds
\begin{equation}
\begin{array}{l}
W_{n+1}^{\epsilon}(g,r_{n+1}) = (k_{\epsilon} \ominus \tilde{W}_{n}^{\epsilon}(\cdot,r_n))(g) + \epsilon \\
 = \left\{
 \begin{array}{ll}
 \sum \limits_{m=0}^{n} (m+1) \epsilon \, 1_{[r_{m}, r_{m+1}]}(d(g,e)) &\textrm{if }d(g,e)\leq r_{n+1}=(n+1)\epsilon \\
 \infty &\textrm{if }d(g,e)> r_{n+1},
 \end{array}
 \right.
 \end{array}
\end{equation}
where $1_{[r_{m}, r_{m+1}]}$ denotes the indicator function on set $[r_m,r_{m+1}]$.
\end{corollary}

\textbf{Proof }The first part follows by Lemma~\ref{lemma:imp} for $r=r_{n+1}$ (i.e. $r_{new}=\epsilon$). The second part follows by induction.
Recall from Lemma~\ref{lemma:imp} that $W_{1}^{\epsilon}(g,r)=k_r(g)+r$.
Now application of (\ref{mainidea}) for $n=1$ yields
{\small
\begin{equation} \label{init}
\begin{array}{l}
W_{2}^{\epsilon}(g,r_2)= (k_{\epsilon}\ominus \tilde{W}_{1}^{\epsilon}(\cdot,r_1))(g) + \epsilon=
\epsilon + \inf \limits_{h \in \overline{B_{g,\epsilon}}} \left\{
\begin{array}{ll}
k_{\epsilon}(h) +\epsilon & \textrm{if }h \neq e \\
0 & \textrm{if }h=e,
\end{array}
\right.
=
\\
\epsilon+ \left\{
\begin{array}{ll}
0 & \textrm{if }d(g,e)\leq \epsilon \\
\epsilon & \textrm{if }\epsilon <d(g,e)\leq 2\epsilon \\
\infty & \textrm{else}
\end{array}
\right.
 =
\left\{
\begin{array}{ll}
\sum \limits_{m=0}^{1} (m+1)\epsilon \; 1_{[r_m,r_{m+1}]}(d(g,e)) & \textrm{if }d(g,e)\leq r_{2} \\
\infty & \textrm{else},
\end{array}
\right.
\end{array}
\end{equation}
}
with $\overline{B_{g,\epsilon}}=\{h \in SE(2)\;|\; d(g,h) \leq \epsilon\}$.
This can intuitively be seen from the geometric meaning of an erosion $\tilde{W}_{1}^{\epsilon} \mapsto
k_{\epsilon}\ominus \tilde{W}_{1}^{\epsilon}$ where one drops cylinders from below on the graph of
$\tilde{W}_{1}^{\epsilon}(\cdot,r_n)$ and considering the new hull where cilinders get stuck.
Eq.~\!(\ref{init}) can also be seen directly from the definition of $k_{\epsilon}$. Let us verify each case separately:
\begin{itemize}
\item If $d(g,e)>2\epsilon$ we
have that the value must be infinite, since suppose it were finite then by the definition of the morphological kernel
$k_{\epsilon}$ we would need to have 
that $d(g,e) \leq d(g,h)+d(h,e) \leq 2\epsilon$ yielding contradiction.
\item If $d(g,e)\leq \epsilon$, then in the erosion-minimization we can take $h=e$ and we obtain $\epsilon+0$.
\item If $\epsilon < d(g,e) \leq 2\epsilon$, then in the erosion-minimization we cannot take $h=e$, but for allowed choices
we obtain $k_{\epsilon}(e)=0$ and
$\epsilon +\epsilon$ as output.
\end{itemize}
Similarly we have by inserting induction hypothesis for $n$ and recursion (\ref{mainidea}) we have
\[
\begin{array}{l}
W_{n+2}^{\epsilon}(g,r_{n+2})= (k_{\epsilon}\ominus \tilde{W}^{\epsilon}_{n+1}(\cdot,r_{n+1}))(g) + \epsilon
=\epsilon+ \sum \limits_{m=0}^{n+1} (m+1)\epsilon \; 1_{[r_{m+1},r_{m+2}]}(d(g,e)) \\
=\epsilon+ \sum \limits_{m'=1}^{n+2} m'\epsilon \; 1_{[r_{m'},r_{m'+1}]}(d(g,e))
= \sum \limits_{m'=0}^{n+1} (m'+1)\epsilon \; 1_{[r_{m'},r_{m'+1}]}(d(g,e)),
\end{array}
\]
for $d(g,e)\leq r_{n+2}$. Here we applied $m'=m+1$ so that the result follows for $n+1$. $\hfill \Box$
\begin{theorem}
Let $g \in SE(2)$ be given.
We have the following limit
\[
\lim \limits_{\epsilon \to 0}
\lim \limits_{n \to \infty} W_{n+1}^{\epsilon}(g,(n+1)\epsilon)
=d(g,e).
\]
\end{theorem}
\ \\
\textbf{Proof }Application of Corollary~\ref{corr:x} gives
\[
\begin{array}{ll}
\lim \limits_{n \to \infty} W_{n+1}^{\epsilon}(g,(n+1)\epsilon) &=\sum \limits_{k=0}^{\infty} (k+1) \epsilon \; 1_{[r_{k},r_{k+1}]}(d(g,e))=\sum \limits_{k=0}^{N^*(g,\epsilon)} (k+1) \epsilon \; 1_{[r_{k},r_{k+1}]}(d(g,e)),
\end{array}
\]
with $N^{*}(g,\epsilon)=\lceil \frac{d(g,e)}{\epsilon}\rceil$, i.e. the largest integer $\geq \frac{d(g,e)}{\epsilon} \in \R^+$. As a result we have
\[
\begin{array}{l}
\lim \limits_{\epsilon \downarrow 0} \lim \limits_{n \to \infty} W_{n+1}^{\epsilon}(g,r_{n+1})=
\lim \limits_{\epsilon \downarrow 0} W_{N^{*}(g,\epsilon)+1}^{\varepsilon}(g, (n+1)\epsilon)
=
\lim \limits_{\epsilon \downarrow 0} \sum \limits_{k=0}^{N^{*}(g,\epsilon)} (k+1)\epsilon \; 1_{[r_{k},r_{k+1}]}(d(g,e))= d(g,e)
\end{array}
\]
where the size of the steps in the staircase towards $d(g,e)$ vanishes as $\epsilon \to 0$. Recall Figure~\ref{fig:hot}. Together with Theorem~\ref{th:3} yielding $W^{\infty}(g)=d(g,e)$, the result follows. $\hfill \Box$


\end{document}